\documentclass[11pt,twoside]{article}
\newcommand{\ifims}[2]{#1}   
\newcommand{\ifAMS}[2]{#1}   
\newcommand{\ifau}[3]{#2}  

\newcommand{\ifbook}[2]{#1}   

\ifbook{

  }
  {

  }

\def\thetitle{Critical dimension in profile semiparametric estimation}
\def\theruntitle{critical dimension in profile semiparametric estimation}

\def\theabstract{
This paper revisits the classical inference results for profile quasi maximum likelihood 
estimators (profile MLE) in the semiparametric estimation problem.
We mainly focus on two prominent theorems: 
the Wilks phenomenon and Fisher expansion for the profile MLE are stated in a new fashion 
allowing finite samples and model misspecification. 
The method of study is also essentially different from the usual analysis 
of the semiparametric problem based on the notion of the hardest parametric submodel.
Instead we derive finite sample deviation bounds for the linear approximation error for the gradient of the loglikelihood.
This novel approach particularly 
allows to address the important issue of the effective target and 
nuisance dimension.
The obtained nonasymptotic results are surprisingly sharp and yield the classical 
asymptotic statements including the asymptotic normality and efficiency of the profile 
MLE. 
The general results are specified to the important special cases of an i.i.d. sample and the analysis is exemplified with a single index model.} 

\def\kwdp{62F10}
\def\kwds{62J12,62F25,62H12}

\def\thekeywords{maximum likelihood, local quadratic bracketing, spread,
local concentration}

\def\authorb{Vladimir Spokoiny}
\def\runauthorb{spokoiny, v.}
\def\addressb{
    Weierstrass Institute and HU Berlin, \\ Moscow Institute of 
    Physics and Technology \\
    Mohrenstr. 39, \\
    10117 Berlin, Germany}
\def\emailb{spokoiny@wias-berlin.de}
\def\affiliationb{Weierstrass-Institute, Humboldt University Berlin, and
Moscow Institute of Physics and Technology}
\def\thanksb
{The author is partially supported by
Laboratory for Structural Methods of Data Analysis in Predictive Modeling, MIPT, 
RF government grant, ag. 11.G34.31.0073.
Financial support by the German Research Foundation (DFG) 
through the CRC 649 ``Economic Risk'' is gratefully acknowledged. 
}

\def\authora{Andreas Andresen}
\def\runauthora{andresen, a.}
\def\addressa{
    \\ Weierstrass-Institute, \\
    Mohrenstr. 39, \\
    10117 Berlin, Germany     
    }
\def\emaila{andresen@wias-berlin.de}
\def\affiliationa{Weierstrass-Institute}
\def\thanksa{The author is supported by Research Units 1735 
"Structural Inference in Statistics: Adaptation and Efficiency"
}

\date{}

\usepackage{amsmath,amssymb,amsthm}
\usepackage{natbib}
\usepackage{epsfig,graphicx}
\usepackage{comment}
\usepackage{color}
\usepackage{srcltx}
\usepackage[mathscr]{eucal}
\usepackage[math]{easyeqn}
\usepackage{etoolbox}

\ifims{
\textheight=23cm
\textwidth=14.8cm
\topmargin=0pt
\oddsidemargin=1.0cm
\evensidemargin=1.0cm
\linespread{1.3}
\renewenvironment{abstract}
    {\centerline{\textbf{Abstract}}\bigskip
      \begin{center}
       \begin{minipage}{11cm}
        \begin{small}
    }
    {   \end{small}
       \end{minipage}
      \end{center}
     \bigskip
    }

}{ 
}

\numberwithin{equation}{section}
\numberwithin{figure}{section}
\newcounter{example}[section]
\numberwithin{example}{section}
\newcounter{remark}[section]
\numberwithin{remark}{section}
\newtheorem{theorem}{Theorem}[section]
\newtheorem{proposition}[theorem]{Proposition}
\newtheorem{lemma}[theorem]{Lemma}
\newtheorem{corollary}[theorem]{Corollary}

\newtheorem{exmp}[example]{Example}
\newtheorem{rmrk}[remark]{Remark}
\newenvironment{example}{\begin{exmp}\rm}{\end{exmp}}
\newenvironment{remark}{\begin{rmrk}\rm}{\end{rmrk}}

\begin{document}
\thispagestyle{empty}
\ifims{
\title{\thetitle}
\ifau{ 
  \author{
    \authora
    \ifdef{\thanksa}{\thanks{\thanksa}}{}
    \\[5.pt]
    \addressa \\
    \texttt{ \emaila}
  }
}
{  
  \author{
    \authora
    \ifdef{\thanksa}{\thanks{\thanksa}}{}
    \\[5.pt]
    \addressa \\
    \texttt{ \emaila}
    \and
    \authorb
    \ifdef{\thanksb}{\thanks{\thanksb}}{}
    \\[5.pt]
    \addressb \\
    \texttt{ \emailb}
  }
}
{   
  \author{
    \authora
    \ifdef{\thanksa}{\thanks{\thanksa}}{}
    \\[5.pt]
    \addressa \\
    \texttt{ \emaila}
    \and
    \authorb
    \ifdef{\thanksb}{\thanks{\thanksb}}{}
    \\[5.pt]
    \addressb \\
    \texttt{ \emailb}
    \and
    \authorc
    \ifdef{\thanksc}{\thanks{\thanksc}}{}
    \\[5.pt]
    \addressc \\
    \texttt{ \emailc}
  }
}

\maketitle
\pagestyle{myheadings}
\markboth
 {\hfill \textsc{ \small \theruntitle} \hfill}
 {\hfill
 \textsc{ \small
 \ifau{\runauthora}
      {\runauthora and \runauthorb}
      {\runauthora, \runauthorb, and \runauthorc}
 }
 \hfill}
\begin{abstract}
\theabstract
\end{abstract}

\ifAMS
    {\par\noindent\emph{AMS 2000 Subject Classification:} Primary \kwdp. Secondary \kwds}
    {\par\noindent\emph{JEL codes}: \kwdp}

\par\noindent\emph{Keywords}: \thekeywords
} 
{ 
\begin{frontmatter}
\title{\thetitle}


\runtitle{\theruntitle}

\ifau{ 
\begin{aug}
    \author{\authora\ead[label=e1]{\emaila}}
    \address{\addressa \\
     \printead{e1}}
\end{aug}

 \runauthor{\runauthora}
\affiliation{\affiliationa} }
{ 
\begin{aug}
    \author{\authora\ead[label=e1]{\emaila}\thanksref{t21}}
    \and
    \author{\authorb\ead[label=e2]{\emailb}\thanksref{t22}}
    
    \address{\addressa \\
     \printead{e1}}
    \address{\addressb \\
     \printead{e2}}
    \thankstext{t21}{\thanksa}
    \thankstext{t22}{\thanksb}
    \affiliation{\affiliationa, \affiliationb} 
    \runauthor{\runauthora and \runauthorb}
\end{aug}
} 
{ 
\begin{aug}
    \author{\authora\ead[label=e1]{\emaila}\thanksref{t21}}
    \and
    \author{\authorb\ead[label=e2]{\emailb}\thanksref{t22}}
    \and
    \author{\authorc\ead[label=e3]{\emailc}\thanksref{t23}}
    
    \address{\addressa \\
     \printead{e1}}
    \address{\addressb \\
     \printead{e2}}
    \address{\addressc \\
     \printead{e3}}
    \thankstext{t21}{\thanksa}
    \thankstext{t22}{\thanksb}
    \thankstext{t23}{\thanksc}
    \affiliation{\affiliationa, \affiliationb, \affiliationc} 
    \runauthor{\runauthora, \runauthorb, and \runauthorc}
\end{aug}}

\begin{abstract}
\theabstract
\end{abstract}

\begin{keyword}[class=AMS]
\kwd[Primary ]{\kwdp}
\kwd[; secondary ]{\kwds}
\end{keyword}

\begin{keyword}
\kwd{\thekeywords}
\end{keyword}

\end{frontmatter}
} 

\def\ND{\cc{N}}
\def\Bernoulli{\mathrm{Bernoulli}}
\def\Vola{\mathrm{Vola}}
\def\Poisson{\mathrm{Poisson}}
\def\ag{\mathrm{ag}}
\def\glob{\operatorname{glob}}
\def\blk{\operatorname{block}}
\def\lin{\operatorname{lin}}
\def\cond{\, \big| \,}

\def\rdl{\epsilon}
\def\rd{\bb{\rdl}}
\def\rddelta{\delta}
\def\rdomega{\varrho}
\def\rddeltab{\rddelta^{*}}
\def\rhorb{\rhor^{*}}

\def\wv{\bb{w}}
\def\varthetav{\bb{\vartheta}}
\def\Lr{\breve{L}}
\def\zetavr{\breve{\zetav}}
\def\etavr{\breve{\etav}}
\def\xivr{\breve{\xiv}}

\def\rdb{\rd}
\def\rdm{\underline{\rdb}}

\def\taub{\tau_{\rdb}}
\def\taum{\tau_{\rdm}}
\def\kappab{\kappa_{\rd}}
\def\deltab{\delta_{\rd}}

\def\taubGP{\tau_{\rdb,\GP}}
\def\taumGP{\tau_{\rdm,\GP}}
\def\kappabGP{\kappa_{\rd,\GP}}
\def\deltabGP{\delta_{\rd,\GP}}
\def\nubm{\nu_{\rd}}
\def\uub{u_{\rd}}
\def\uubGP{u_{\rd,\GP}}
\def\nubmGP{\nu_{\rd, G}}

\def\rG{\rd,\GP}

\def\LinSp{\mathrm{L}}
\def\Id{I\!\!\!I}
\def\Ind{\operatorname{1}\hspace{-4.3pt}\operatorname{I}}

\def\BG{\mathcal{R}}
\def\bg{r}
\def\fmup{\phi}
\def\rg{r}
\def\uc{u_{c}}
\def\muc{\mu_{c}}
\def\mud{\mu_{0}}
\def\xxd{\xx_{0}}
\def\yyd{\yy_{0}}
\def\gmd{\gm_{0}}

\def\ms{m^{*}}
\def\Inv{A}
\def\InvT{\Inv^{\T}}
\def\Invt{\tilde{\Inv}}

\def\ssize{N}
\def\nsize{{n}}

\def\rhor{\omega}

\def\LT{L}
\def\LGP{\LT_{\GP}}
\def\La{\mathbb{L}}
\def\Lab{\La_{\rdb}}
\def\Lam{\La_{\rdm}}

\def\DP{D}
\def\DPc{\DP_{0}}
\def\DPb{\DP_{\rdb}}
\def\DPm{\DP_{\rdm}}

\def\LabGP{\La_{\rdb,\GP}}
\def\LamGP{\La_{\rdm,\GP}}

\def\DPbGP{\DP_{\rdb,\GP}}
\def\DPmGP{\DP_{\rdm,\GP}}
\def\riskbGP{\riskt_{\rdb,\GP}}

\def\gmi{\mathtt{b}}
\def\gmiid{\mathtt{g}_{1}}
\def\kullbi{\Bbbk}
\def\Thetasi{\Theta_{\loc}}
\def\rri{\mathtt{u}}
\def\rris{\rri_{0}}

\def\Ipc{\bb{\mathrm{f}}}
\def\IF{\Bbb{F}}
\def\IFc{\IF_{0}}
\def\IFb{\IF_{\rdb}}
\def\IFm{\IF_{\rdm}}

\def\DF{\cc{D}}
\def\DFc{\DF_{0}}
\def\DFb{\DF_{\rdb}}
\def\DFm{\breve{\DF}_{\rd}}
\def\DFm{\DF_{\rdm}}

\def\DPr{\breve{\DP}}
\def\VF{\cc{V}}
\def\VFc{\VF_{0}}

\def\HHc{\HH_{0}}
\def\HHb{\HH_{\rd}}
\def\HHm{\HH_{\rdm}}

\def\xib{\xi^{*}}
\def\xivb{\xiv_{\rdb}}
\def\xivm{\xiv_{\rdm}}
\def\CAm{\underline{\CA}}
\def\CAb{\CA}

\def\penr{\operatorname{pen}}
\def\pen{\mathfrak{t}}
\def\PEN{\operatorname{PEN}}
\def\RSS{\operatorname{RSS}}
\def\med{\operatorname{med}}

\def\ex{\mathrm{e}}
\def\entrl{\mathbb{Q}}
\def\entrlb{\entrl}
\def\entr{\entrl}

\def\kullb{\cc{K}} 
\def\kullbc{\kullb^{c}}

\def\gm{\mathtt{g}}
\def\gmc{\gm_{c}}
\def\gmb{\gm}
\def\gmbm{\gmb_{1}}

\def\yy{\mathtt{y}}
\def\yyc{\yy_{c}}
\def\xx{\mathtt{x}}
\def\xxc{\xx_{c}}
\def\tc{t_{c}}

\def\alp{\alpha}
\def\alpn{\rho}
\def\gmu{\mathfrak{a}}

\def\losst{\varrho}
\def\loss{\wp}
\def\lossp{u}
\def\closs{g}

\def\riskt{\cc{R}}
\def\emprisk{\ell}
\def\bias{b}
\def\bern{q}

\def\TT{\nsize}

\def\Pone{P}
\def\Pf{\P_{f(\cdot)}}
\def\Ef{\E_{f(\cdot)}}
\def\Ps{\P_{\thetas}}
\def\Es{\E_{\thetas}}
\def\Pu{\P_{\upsilons}}
\def\Eu{\E_{\upsilons}}

\def\Pvs{\P_{\thetavs}}
\def\Evs{\E_{\thetavs}}

\def\UPd{w}
\def\nunup{\nu_{1}}
\def\rru{\rr_{1}}
\def\rups{\rr_{0}}
\def\rupsb{\rups^{*}}
\def\rrf{\rr^{\flat}}

\def\smooths{\mathbb{S}}
\def\smooth{\smooths_{1}}

\def\elli{\bar{\ell}}

\def\K{K}

\def\Psir{\breve{\Psi}}

\def\af{a}
\def\afs{\af^{*}}

\def\kapla{\varkappa}

\newcommand{\mlew}[1]{\tilde{\thetav}_{#1}}
\newcommand{\mlea}[1]{\hat{\thetav}_{#1}}
\newcommand{\mluw}[1]{\tilde{\theta}_{#1}}
\newcommand{\mlua}[1]{\hat{\theta}_{#1}}
\newcommand{\penm}[1]{\boldsymbol{m}_{#1}}

\def\Pdom{\mu_{0}}
\def\PDOM{\bb{\mu}_{0}}
\def\EDOM{\E_{0}}

\def\mk{m}
\def\Mk{\cc{M}}
\def\SV{\cc{S}}

\def\Cs{E}
\def\Csd{\Cs^{\circ}}
\def\Ca{A}
\def\CS{\cc{E}}
\def\CA{\cc{A}}
\def\CAb{\CA_{\rd}}
\def\CAC{\CA_{\CoFu}}

\def\Ccb{m_{\rdb}}
\def\Ccm{m_{\rdm}}
\def\CcbGP{m_{\rdb,\GP}}
\def\CcmGP{m_{\rdm,\GP}}

\def\etas{\eta^{*}}

\def\omegav{\bb{\phi}}
\def\omegavs{\omegav^{*}}
\def\omegavc{\omegav'}

\def\nuvs{\nuv^{*}}
\def\nuvc{\nuv'}

\def\nunu{\nu_{0}}
\def\numu{\nu_{1}}
\def\nupi{\nu^{+}}
\def\nubu{\beta}

\def\nus{\nu}
\def\nusb{\nus}
\def\nusr{\nus^{\bracketing}}
\def\Nusb{\mathbb{N}}
\def\Nusr{\mathbb{N}^{\diamond}}

\def\dist{d}
\def\distd{\mathfrak{a}}

\def\hatk{\kappa}
\def\ko{k^{\circ}}

\def\qqq{\mathfrak{q}}
\def\ppp{{s}}
\def\Cqq{C(\qqq)}
\def\Cqqb{C^{\diamond}(\qqq)}
\def\Crho{C(\mrho)}
\def\Cqqm{\log(4)}
\def\Cqpr{(\qqq \rrp + \dimp / 2)}

\def\Cdima{\mathfrak{e}_{0}}
\def\Cdimb{\mathfrak{e}_{1}}
\def\cdima{\mathfrak{c}_{0}}
\def\cdimb{\mathfrak{c}_{1}}
\def\cdim{\mathfrak{c}}

\def\rdomega{\varrho}
\def\deltaD{\delta}
\def\alphai{\alpha_{1}}
\def\alphaii{\alpha_{2}}
\def\alphaiii{\alpha_{3}}
\def\alphaiv{\alpha_{4}}

\def\err{\diamondsuit}
\def\errbm{\bar{\err}_{\rdomega}}
\def\errm{\err_{\rdm}}
\def\errb{\err_{\rdb}}

\def\errbGP{\err_{\rdomega,\GP}}
\def\errmGP{\err_{\rdm,\GP}}
\def\errbmGP{\bar{\err}_{\rd,\GP}}

\def\errs{\err_{\rdomega}^{*}}
\def\deltas{\alpha}

\def\xivbGP{\xiv_{\rdb,\GP}}
\def\xivmGP{\xiv_{\rdm,\GP}}

\def\SP{S}
\def\GP{G}
\def\GPt{\GP_{0}}
\def\GPn{\GP_{1}}
\def\gp{g}
\def\gs{s}

\def\errbGP{\err_{\rdb,\GP}}
\def\errmGP{\err_{\rdm,\GP}}
\def\errpmGP{\err_{\GP}^{\pm}}

\def\LCS{\cc{C}}

\def\DPGP{\DP_{\GP}}
\def\thetavsGP{\thetavs_{\GP}}

\def\LL{\cc{L}}
\def\LLb{\LL^{*}}
\def\LLh{\cc{L}}

\def\YY{\cc{Y}}
\def\LP{L^{\circ}}

\def\modcnrd{\mathfrak{A}}

\def\pens{\pi}
\def\pnn{\mathfrak{g}}
\def\pnnd{\mathfrak{u}}
\def\pnndGP{\pnnd_{\GP}}

\def\confpr{\mathfrak{c}}
\def\confprb{\confpr^{*}}

\def\pn{\pens^{*}}
\def\penInt{\mathfrak{D}}
\def\penH{\mathbb{H}}
\def\pmu{\mathfrak{u}}
\def\Closs{\cc{R}}

\def\dimp{p}
\def\riskb{\riskt_{\rdb}}
\def\dimpp{\dimp+1}
\def\BB{I\!\!B}
\def\vA{\mathtt{v}}

\def\deficiency{\Delta}
\def\spread{\Delta}
\def\dimtotal{\dimp^{*}}

\def\thetav{\bb{\theta}}
\def\thetavs{\thetav^{*}}
\def\thetavc{\thetav'}
\def\thetavd{\thetav^{\circ}}
\def\thetavdc{\thetav^{\sharp}}
\def\dthetavs{\thetav,\thetavs}

\def\thetas{\theta^{*}}
\def\thetac{\theta'}
\def\thetad{\theta^{\circ}}
\def\thetab{\theta^{\dag}}
\def\thetavb{\thetav^{\dag}}

\def\vtheta{\vartheta}
\def\vthetav{\bb{\vtheta}}
\def\prior{\Pi}

\def\Gam{\Xi}
\def\Gam{\mathcal{S}}
\def\RG{R}
\def\Psu{\Upsilon}
\def\Phim{\breve{\Phi}}

\def\Proj{P}

\def\gammavs{\gammav^{*}}
\def\gammavd{\gammav^{\circ}}
\def\etavs{\etav^{*}}
\def\etavd{\etav^{\circ}}
\def\etavc{\etav'}

\def\taus{\tau_{0}}
\def\taup{\lceil \tau \rceil}

\def\sigmas{{\sigma^{*}}}
\def\Sigmas{\Sigma_{0}}

\def\upsilonc{\upsilon'}
\def\upsilond{\upsilon^{\circ}}
\def\upsilonp{{\upsilon}^{*}}
\def\upsilonm{{\upsilon}_{*}}
\def\upsilonvs{\upsilonv^{*}}
\def\upsilons{\upsilon^{*}}
\def\upsilonb{\bar{\upsilon}}
\def\upsilonvd{\upsilonv^{\circ}}

\def\ups{\bb{\upsilon}}
\def\upss{\ups_{0}}
\def\upsc{\ups^{\prime}}
\def\upsd{\ups^{\circ}}
\def\upsdc{\ups^{\sharp}}
\def\upsdu{\ups^{\flat}}

\def\Ups{\varUpsilon}
\def\Upsd{\Ups^{\circ}}
\def\Upss{\Ups_{\circ}}
\def\UpsP{\Ups^{c}}

\def\Thetas{\Theta_{0}}
\def\ThetasGP{\Theta_{0,\GP}}
\def\varthetav{\bb{\vartheta}}

\def\glink{g}

\def\fvs{\fv}
\def\fs{f}
\def\fb{\fv^{\dag}}

\def\uc{\uv'}
\def\ud{\uv^{\circ}}
\def\uvs{\uv^{*}}
\def\us{u^{*}}
\def\vs{v^{*}}

\def\reps{\epsilon}
\def\eps{\epsilon}

\def\repsc{\reps_{0}}
\def\repsb{\reps^{*}}
\def\repsg{g}

\def\lu{\delta}
\def\lub{\bar{\lu}}

\def\Uu{U}
\def\UU{\cc{Y}}
\def\UUM{\cc{M}}
\def\UP{\cc{U}}
\def\up{\mathfrak{u}}

\def\VP{V}
\def\VPc{\VP_{0}}
\def\VPV{\cc{U}}
\def\VPVc{\cc{\VPV}_{0}}
\def\VPGP{\VP_{\GP}}
\def\VPSP{\VP_{\SP}}

\def\VV{H}
\def\GV{\cc{G}}
\def\GVS{S}

\def\VVb{\VV^{*}}
\def\VVc{\VV_{0}}
\def\vv{\bb{h}}
\def\vva{g}
\def\vp{\mathbf{v}}
\def\vpc{\vp_{0}}
\def\VVca{\VV}
\def\Vtt{H}

\def\DG{D}

\def\Vn{V_{0}}
\def\vn{v_{0}}

\def\norm{\mathfrak{c}}
\def\normc{\delta}
\def\norma{c}

\def\egridd{\cc{E}_{\delta}}
\def\penb{\varkappa}

\def\dotzeta{\dot{\zeta}}
\def\mes{\pi}
\def\mesl{\Lambda}
\def\cprr{F}

\def\lambdam{\gm_{1}}
\def\lambdaB{{\lambda}^{*}}
\def\lambdac{{\lambda'}}

\def\cla{{b}}
\def\fis{\mathfrak{a}}
\def\fiss{\fis_{1}}

\def\Vd{{V}}
\def\vd{\bar{v}}

\def\klim{k^{\circ}}
\def\midm{\mid \!}

\def\Ldrift{M}
\def\ldrift{m}
\def\mY{b}
\def\Lvar{D}
\def\lvar{\sigma}

\def\Mubcu{\Upsilon}
\def\Dthetav{\bb{u}}

\def\B{\cc{B}}
\def\BD{\B^{\circ}}
\def\BU{B}
\def\BI{\B^{*}}

\def\mub{\mu^{*}}
\def\mubc{\mu}
\def\mubcb{\mubc^{*}}
\def\Mubc{\mathbb{M}}
\def\Mubcb{\mathrm{M}}

\def\zzc{\zz_{c}}
\def\ww{w}
\def\wwc{\ww_{c}}

\def\norms{\circ} 
\def\rs{\rr_{\norms}}
\def\yys{\yy_{\norms}}
\def\xxs{\xx_{\norms}}
\def\zzs{\zz_{\norms}}
\def\uu{\mathtt{u}}
\def\uus{\uu_{\norms}}
\def\mus{\mu_{\norms}}
\def\gms{\gm_{\norms}}
\def\wws{\ww_{\circ}}

\def\srho{s}
\def\mrho{\varrho}

\def\Lmgf{\mathfrak{M}}
\def\Lmgfb{\Lmgf^{*}}

\def\lmgf{\mathfrak{m}}
\def\lmgfb{\lmgf^{*}}

\def\Expzeta{\mathfrak{N}}
\def\expzeta{\mathfrak{s}}

\def\rr{\mathtt{r}}
\def\rrb{\rr^{*}}
\def\rru{\rr_{\circ}}
\def\rrc{\rr'}
\def\rs{r_{*}}

\def\zz{\mathfrak{z}}
\def\zzb{\tilde{\zz}}
\def\tt{\mathfrak{t}}
\def\zb{z_{\rd}}
\def\zzg{\zz_{1}}
\def\zzQ{\zz_{0}}
\def\zzq{\zz}

\def\Cr{\mathfrak{c}}
\def\Crp{\mathfrak{C}}
\def\Crl{\mathfrak{r}}
\def\Crlp{\mathfrak{R}}
\def\Crlq{\cc{T}}
\def\Crlmu{\cc{M}}

\def\zetah{\zeta_{h}}
\def\GG{G}
\def\HH{H}
\def\pG{p}
\def\pH{q}
\def\hh{H^{*}}

\def\mubch{\mubc_{1}}
\def\rhoh{\rho_{1}}
\def\CoFuh{\CoFu_{1}}
\def\dimh{p_{1}}
\def\VPh{\VP_{1}}
\def\VPt{\VP_{0}}

\def\LLh{L_{1}}
\def\pnndh{\pnnd_{1}}

\def\LCS{C}
\def\Ac{A_{0}}
\def\Ab{A_{\rd}}
\def\DPrb{\DPr_{\rdb}}
\def\DPrm{\DPr_{\rdm}}
\def\Cb{\cc{C}_{\rdb}}
\def\Ub{\cc{U}_{\rdb}}
\def\zetavrb{\zetavr_{\rd}}
\def\xivrb{\breve{\xiv}_{\rd}}
\def\VPrb{\breve{\VP}_{\rdb}}
\def\Larb{\breve{\La}_{\rdb}}
\def\Larm{\breve{\La}_{\rdm}}

\def\deltav{\bb{\delta}}

\def\score{\nabla}
\def\scorer{\breve{\nabla}}

\def\LCS{C}
\def\Ac{A_{0}}
\def\Bc{B_{0}}
\def\AF{A}
\def\Ab{A_{\rdb}}
\def\Am{A_{\rdm}}
\def\DPrc{\DPr_{0}}
\def\DPrb{\DPr_{\rdb}}
\def\DPrm{\DPr_{\rdm}}
\def\Cb{\cc{C}_{\rdb}}
\def\Cm{\cc{C}_{\rdm}}
\def\Ub{\cc{U}_{\rdb}}
\def\deltav{\bb{\delta}}
\def\nuv{\bb{\nu}}
\def\xivrb{\breve{\xiv}_{\rd}}
\def\VPrb{\breve{\VP}_{\rdb}}
\def\Larb{\breve{\La}_{\rdb}}
\def\Lar{\breve{\La}}
\def\Larm{\breve{\La}_{\rdm}}
\def\VH{Q}
\def\VHc{\VH_{0}}
\def\zetavrm{\zetavr_{\rdm}}
\def\N{\mathbb{N}}

\def\Span{\operatorname{span}}
\def\Exc{{\square}}
\def\UUs{U_{\circ}}
\def\errbm{\errb^{*}}
\def\corrDF{\nu}
\def\BBr{\breve{\BB}}
\def\taua{\tau}
\def\AssId{\mathcal{I}}
\def\assId{\iota}
\def\AFD{\cc{A}}

\def\BanX{\cc{X}}
\def\basX{\ev}
\def\apprX{\alpha}
\def\fvs{\fv^{*}}
\def\lkh{\ell}
\def\Bc{B_{0}}
\def\dimn{\dimp_{\nsize}}
\def\betan{\beta_{\nsize}}


\def\xivGP{\xiv_{\GP}}
\def\dimA{\mathtt{p}}
\def\dimAGP{\dimA}
\def\dime{\dimA_{e}}
\def\dimG{\dimA_{\GP}}
\def\dimS{\dimA_{s}}
\def\nubm{\nu_{\rd}}
\def\uub{u_{\rd}}
\def\uubGP{u_{\rd,\GP}}

\def\priorden{\pi}
\def\xivGP{\xiv_{\GP}}
\def\dimAGP{\dimA}
\def\nubm{\nu_{\rd}}
\def\uub{u_{\rd}}
\def\uubGP{u_{\rd,\GP}}

\def\CR{\mathcal{C}}
\def\CRb{\CR_{\rdb}}
\def\vthetavb{\bar{\vthetav}}
\def\Covpost{\mathfrak{S}}

\def\Db{\DP_{+}}
\def\Dm{\DP_{-}}
\def\uvb{\uv_{+}}
\def\uvm{\uv_{-}}
\def\uud{\omega}
\def\taub{\delta}
\def\Lip{L}
\def\Xb{X_{+}}
\def\Xm{X_{-}}
\def\deltam{\delta_{-}}
\def\betauv{\delta}
\def\betab{\betauv_{1}}
\def\betaf{\betauv_{2}}
\def\upsv{\bb{\varkappa}}
\def\upsvb{\bar{\upsv}}
\def\rhob{\varrho}
\def\alpb{\alp_{1}}
\def\betap{\betauv_{3}}
\def\Ec{\E^{\circ}}
\def\ff{f}
\def\fpos{g}
\def\fneg{h}
\def\alpb{\alp_{+}}
\def\alpm{\alp_{-}}

\def\kappak{\kappa}
\def\kappas{\kappak^{*}}
\def\Kappak{\cc{K}}
\def\DPk{\DP_{\kappak}}
\def\VPk{\VP_{\kappak}}

\def\ts{s}
\def\tsv{\bb{\ts}}
\def\mm{\kappa}
\def\mmc{\mm'}
\def\mmd{\mm^{\circ}}
\def\mmo{\mm^{*}}
\def\mmmmo{\mm,\mmo}
\def\mmt{\tilde{\mm}}
\def\mma{\hat{\mm}}
\def\pp{z}

\def\LLL{L_{1}}
\def\LLr{L_{0}}
\def\muL{\mu_{1}}
\def\mur{\mu_{0}}

\def\LmgfL{\Lmgf_{1}}
\def\Lmgfr{\Lmgf_{0}}
\def\Lmgfm{\Lmgf_{1}}

\def\Kappa{\cc{K}}
\def\CoFu{\cc{C}}
\def\CoFuc{\CoFu_{0}}
\def\CoFub{\CoFu^{*}}
\def\CoFuL{\CoFu_{1}}
\def\CoFur{\CoFu_{0}}
\def\CAL{\CA_{1}}
\def\CAr{\CA_{0}}
\def\CAzz{\cc{A}}

\def\pnnL{\pnn_{1}}
\def\pnnr{\pnn_{0}}
\def\ttd{\delta}
\def\alphaL{\alpha_{1}}
\def\alphar{\alpha_{0}}
\def\alpharL{\alpha}
\def\rat{\mathfrak{t}}
\def\mquad{\nquad}
\def\zzL{\zz_{1}}
\def\zzr{\zz_{0}}

\def\mmset{\mathcal{I}}
\def\xex{u}
\def\dcm{q}
\def\dc{g}
\def\dcL{\dc_{1}}
\def\dcr{\dc_{0}}
\def\kk{k}

\def\cpen{\tau}

\def\dens{f}
\def\jj{j}
\def\JJ{\cc{J}}
\def\Zphi{Z}
\def\Zphiv{\bb{\Zphi}}

\def\nuu{\mathfrak{u}}
\def\nud{\mathfrak{u}_{0}}
\def\nun{c_{\nuu}}
\def\rhork{\kullb}
\def\GH{\mbox{GH}}
\def\HYP{\mbox{HYP}}
\def\NIG{\mbox{NIG}}
\def\IR{{\rm I\!R}}
\def\taggr{b}
\def\penm{\boldsymbol{m}}
\def\Crlp{\cc{R}}

\def\Mh{M}
\def\Mht{\Mh^{c}}

\def\Mhh{\Mh^{-}}
\def\Mhc{G}
\def\Lh{L_{1}}
\def\Uh{\cc{U}}
\def\wloc{w}
\def\Bias{B}
\def\bias{b}
\def\ExpzetaU{\Expzeta_{1}}
\def\vpci{\vp_{i,0}}
\def\IFci{\IF_{i,0}}

\def\erqb{\Circle_{\rdb}}
\def\erqm{\Circle_{\rdm}}
\def\errqm{\errm^{*}}
\def\errqb{\errb^{*}}
\def\Nsize{N}
\def\VVD{\VV_{1}}
\def\AA{A}
\def\Wloc{W}

\renewcommand{\(}{$\,}
\renewcommand{\)}{\,$}

\def\nquad{\hspace{-1cm}}
\def\eqdef{\stackrel{\operatorname{def}}{=}}
\def\tod{\stackrel{d}{\longrightarrow}}
\def\tow{\stackrel{w}{\longrightarrow}}
\def\toP{\stackrel{\P}{\longrightarrow}}

\newcommand{\cc}[1]{\mathscr{#1}}
\newcommand{\bb}[1]{\boldsymbol{#1}}

\renewcommand{\bar}[1]{\overline{#1}}
\renewcommand{\hat}[1]{\widehat{#1}}
\renewcommand{\tilde}[1]{\widetilde{#1}}

\renewcommand{\Gamma}{\varGamma}
\renewcommand{\Pi}{\varPi}
\renewcommand{\Sigma}{\varSigma}
\renewcommand{\Delta}{\varDelta}
\renewcommand{\Lambda}{\varLambda}
\renewcommand{\Psi}{\varPsi}
\renewcommand{\Phi}{\varPhi}
\renewcommand{\Theta}{\varTheta}
\renewcommand{\Omega}{\varOmega}
\renewcommand{\Xi}{\varXi}
\renewcommand{\Upsilon}{\varUpsilon}
\def\nn{\nonumber \\}

\def\suml{\sum\limits}
\def\supl{\sup\limits}
\def\maxl{\max\limits}
\def\infl{\inf\limits}
\def\intl{\int\limits}
\def\liml{\lim\limits}
\def\Cov{\operatorname{Cov}}
\def\Var{\operatorname{Var}}
\def\arginf{\operatornamewithlimits{arginf}}
\def\argsup{\operatornamewithlimits{argsup}}
\def\argmax{\operatornamewithlimits{argmax}}
\def\argmin{\operatornamewithlimits{argmin}}
\def\val{\operatorname{val}}

\def\D{\boldsymbol{D}}
\def\dd{\operatorname{d}}
\def\tr{\operatorname{tr}}
\def\I{I\!\!I}
\def\R{I\!\!R}
\def\E{I\!\!E}
\def\P{I\!\!P}
\def\X{\mathfrak{X}}
\def\kappa{\varkappa}
\def\Const{\mathrm{Const.} \,}
\def\cdt{\boldsymbol{\cdot}}
\def\tm{\!\times\!}
\def\T{\top}
\def\diag{\operatorname{diag}}
\def\diam{\operatorname{diam}}
\def\rank{\operatorname{rank}}
\def\loc{\operatorname{loc}}

\def\av{\bb{a}}
\def\bv{\bb{b}}
\def\cv{\bb{c}}
\def\dv{\bb{d}}
\def\ev{\bb{e}}
\def\fv{\bb{f}}
\def\gv{\bb{g}}
\def\hv{\bb{h}}
\def\iv{\bb{i}}
\def\jv{\bb{j}}
\def\kv{\bb{k}}
\def\lv{\bb{l}}
\def\mv{\bb{m}}
\def\nv{\bb{n}}
\def\ov{\bb{o}}
\def\pv{\bb{p}}
\def\qv{\bb{q}}
\def\rv{\bb{r}}
\def\sv{\bb{s}}
\def\tv{\bb{t}}
\def\uv{\bb{u}}
\def\vv{\bb{v}}
\def\wv{\bb{w}}
\def\xv{\bb{x}}
\def\yv{\bb{y}}
\def\zv{\bb{z}}

\def\Cv{\bb{C}}
\def\Gv{\bb{G}}
\def\Mv{\bb{M}}
\def\Sv{\bb{S}}
\def\Uv{\bb{U}}
\def\Xv{\bb{X}}
\def\Yv{\bb{Y}}
\def\Zv{\bb{Z}}

\def\alphav{\bb{\alpha}}
\def\epsv{\bb{\varepsilon}}
\def\etav{\bb{\eta}}
\def\gammav{\bb{\gamma}}
\def\varepsilonv{\bb{\varepsilon}}
\def\phiv{\bb{\phi}}
\def\psiv{\bb{\psi}}
\def\tauv{\bb{\tau}}
\def\upsilonv{\bb{\upsilon}}
\def\xiv{\bb{\xi}}
\def\zetav{\bb{\zeta}}

\def\Psiv{\bb{\Psi}}
\def\CONST{\mathtt{C}}

\def\itemv{\vfill\item}
\newenvironment{myslide}[1]
    {\begin{frame}\frametitle{#1}\vfill}
    {\vfill\end{frame}}

\def\vsp{\vspace{0.05\textheight} \vfill}
\def\summarysign{\resizebox{0.08\textwidth}{0.08\textheight}{\includegraphics{summary}}\,}
\def\nix{}
\def\wpu{$\bullet$}

\def\btri{\vfill{\( \blacktriangleright \) }}
\def\btrir{\vfill{\( \blacktriangleright \) }}

\newcommand{\mygraphics}[3]{\begin{center}
    \resizebox{#1\textwidth}{#2\textheight}{\includegraphics{#3}}
    \end{center}
}

\newcommand{\mybox}[3]{\begin{center}
    \resizebox{#1\textwidth}{#2\textheight}{#3}
    \end{center}
}

\newenvironment{eqnh}
{
    \setbeamercolor{postit}{fg=black,bg=hellgelb} 
    \begin{beamercolorbox}[center,wd=\textwidth]{postit} 
    \begin{eqnarray*}}
    {\end{eqnarray*}\end{beamercolorbox}
}

\def\gps{s}
\def\GK{\cc{G}}
\def\Excgr{\diamondsuit}

\def\dimh{m}
\def\LCS{C}
\def\Ac{A}
\def\Bc{E}
\def\AF{A}
\def\CF{C}
\def\Ab{A_{\rdb}}
\def\Am{A_{\rdm}}
\def\DPc{\DP}
\def\VPc{\VP}
\def\HHc{\HH}
\def\DPrc{\DPr}
\def\DPrp{\DPr_{\dimh}}
\def\DPrb{\DPr_{\rdb}}
\def\DPrm{\DPr_{\rdm}}
\def\Cb{\cc{C}_{\rdb}}
\def\Cm{\cc{C}_{\rdm}}
\def\Ub{\cc{U}_{\rdb}}
\def\xivrb{\breve{\xiv}_{\rd}}
\def\VPrb{\breve{\VP}_{\rdb}}
\def\Larb{\breve{\La}_{\rdb}}
\def\Lar{\breve{\La}}
\def\Larm{\breve{\La}_{\rdm}}

\def\DFc{\DF}
\def\VFc{\VF}

\def\VH{Q}
\def\VHc{\VH}
\def\zetavrm{\zetavr_{\rdm}}

\def\fvh{\bb{\dimh}}
\def\N{\mathbb{N}}
\def\Z{\mathbb{Z}}

\def\iic{\IF}
\def\iif{\breve{\iic}}
\def\DP{{D}}
\def\HH{{H}}
\def\A{{A}}
\def\ifc{\breve{\iic}}

\def\deltar{\delta}

\def\Thetathetav#1{\substack{\\[0.1pt] \upsilonv \in \Ups \\[1pt] \Proj \upsilonv = #1}}
\def\Span{\operatorname{span}}
\def\Exc{{\square}}
\def\UUs{U_{\circ}}
\def\errbm{\errb^{*}}
\def\corrDF{\rho}
\def\BBr{\breve{\BB}}
\def\taua{\tau}
\def\AssId{\mathcal{I}}
\def\AFD{\cc{A}}

\def\BanX{\cc{X}}
\def\basX{\ev}
\def\apprX{\alpha}
\def\fvs{\fv^{*}}
\def\lkh{\ell}
\def\h{\frac{1}{2}}
\def\basis{\ev}
\def\Proj{\Pi_{0}}

\def\Ij{\mathcal{I}}

\def\Mn{M_{\nsize}}
\def\bA{\breve{A}}
\def\cA{\bA_{\dimh}}

\def\Sdr{\cc{S}}
\def\xxn{\xx_{\nsize}}

\def\CONST{\mathtt{C}}
\def\Ij{\mathcal{I}}

\def\etas{\eta^{*}}
\def\zetavs{\zetav^{*}}
\def\zetavc{\zetav'}

\def\omegav{\bb{\phi}}
\def\omegavs{\omegav^{*}}
\def\omegavc{\omegav'}

\def\dimn{\dimp_{\nsize}}
\def\betan{\beta_{\nsize}}

\def\bA{\breve{A}}
\def\cA{\bA_{\dimh}}

\def\corrDF{\rho}
\def\rupf{\rr_{1}}
\def\nuno{\nunu}
\def\gmone{\gm_{1}}
\def\rhorb{\rhor_{1}}

\def\upsilonv{\boldsymbol{\upsilon}}
\def\upsilonvs{\boldsymbol{\upsilon}^{*}}
\def\upsilonvd{\boldsymbol{\upsilon}^\circ}
\def\upsilonvc{\upsilonv'}

\def\dimB{\mathtt{p}_{\BB}}

\section{Introduction}
\label{Chgsemiparam}
\label{Ssemipar}

Many statistical tasks can be viewed as problems of semiparametric estimation when 
the unknown data distribution is described by a high or infinite dimensional 
parameter while the target is of low dimension. 
Typical examples are provided by functional estimation, estimation of a function at 
a point, or simply by estimating a given subvector of the parameter vector.
The classical statistical theory provides a general solution to this problem:
estimate the full parameter vector by the maximum likelihood method and project the 
obtained estimate onto the target subspace. 
This approach is known as \emph{profile maximum likelihood} and it appears to be 
\emph{semiparametrically efficient} under some mild regularity conditions.
We refer to the papers \cite{MurphyvanderVaart,MuVa1999} and the book \cite{Kosorok} 
for a detailed presentation of the modern state of the theory and further 
references. 
The famous Wilks result claims that the likelihood ratio test statistic 
in the semiparametric test problem is nearly chi-square with \( \dimp \) degrees of 
freedom corresponding to the dimension of the target parameter.
Various extensions of this result can be found e.g. in 
\cite{FaZh2001,FaHu2005,BoMa2011}; see also the references therein.

This study revisits the problem of profile semiparametric estimation and 
addresses some new issues. 
The most important difference between our approach and the classical theory is 
a nonasymptotic character of our study.
A finite sample analysis is particularly challenging because most of notions, 
methods and tools in the classical theory are formulated in the asymptotic setup 
with growing sample size. 
Only few finite sample general results are available; see e.g. the recent paper 
\cite{BoMa2011}.
The results of this paper explicitly describes all ``small'' terms 
in the expansion of the log-likelihood.
This helps to carefully treat the question of applicability of the approach 
in different situations. 
A particularly important question concerns the critical dimension of the target \( \dimp \) 
and the full parameter dimension \( \dimtotal \) for which the main results are 
still accurate. 
Another issue addressed in this paper is the model misspecification. 
In many practical problems, it is unrealistic to expect that the model 
assumptions are exactly fulfilled, even if some rich nonparametric models are used.
This means that the true data distribution \( \P \) does not belong to the 
considered parametric family.
Applicability of the general semiparametric theory in such cases is questionable. 
An important feature of the presented approach is that it equally applies 
under a possible model misspecification.


%
Let \( \Yv \) denote the observed random data, and \( \P \) denote the data 
distribution. 
The parametric statistical model assumes that the unknown data distribution \( \P \) 
belongs to a given parametric family \( (\P_{\upsilonv}) \):
\begin{EQA}[c]
    \Yv
    \sim 
    \P = \P_{\upsilonvs} \in (\P_{\upsilonv}, \, \upsilonv \in \Ups),
\end{EQA}
where \( \Ups \) is some high dimensional or even infinite dimensional parameter 
space. 

The maximum likelihood approach in the parametric 
estimation suggests to estimate the whole parameter vector \( \upsilonv \) by 
maximizing the corresponding log-likelihood 
\( \LL(\upsilonv) = \log \frac{d\P_{\upsilonv}}{d\PDOM}(\Yv) \) for some dominating 
measure \( \PDOM \): 
\begin{EQA}[c]
\label{profileest}
    \tilde{\upsilonv} 
    \eqdef 
    \argmax_{\upsilonv \in \Ups} \LL(\upsilonv).
\end{EQA}
Our study admits a model misspecification 
\( \P \notin (\P_{\upsilonv}\, , \upsilonv \in \Ups) \).
Equivalently, one can say that \( \LL(\upsilonv) \) is the 
\emph{quasi log-likelihood function} on \( \Ups \). 
The ``target'' value \( \upsilonvs \) of the parameter \( \upsilonv \) 
can defined by 
\begin{EQA}[c]
\label{eq: definition of full target}
    \upsilonvs 
    = 
    \argmax_{\upsilonv \in \Ups} \E \LL(\upsilonv) .
\end{EQA}
Under model misspecification, 
\( \upsilonvs \) defines the best parametric fit to \( \P \) by the considered 
family. 

In the semiparametric framework, the target of analysis is only a 
low dimensional component \( \thetav \) of the whole parameter \( \upsilonv \). 
This means that the target of estimation is 
\begin{EQA}[c]
\label{thupvs}
    \thetavs
    =
    \Proj \upsilonvs,
\end{EQA}
for some mapping \( \Proj : \Ups \to \R^{\dimp} \), and \( \dimp \in \N \) stands 
for the dimension of the target. Often the vector \( \upsilonv \) is 
represented as \( \upsilonv = (\thetav,\etav) \), where \( \thetav \) is the target 
of analysis while \( \etav \) is the \emph{nuisance parameter}. 
We refer to this situation as \( (\thetav,\etav) \)-setup and 
our presentation follows this setting.

Define
\begin{EQA}[c]
    \Lr(\thetav) 
    \eqdef 
    \max_{\Thetathetav{\thetav}} \LL(\upsilonv) .
\label{Lrthetav}
\end{EQA}   
The \emph{profile maximum likelihood} approach defines the estimator of \( \thetavs \) 
by projecting the obtained MLE \( \tilde{\upsilonv} \) on the target space:
\begin{EQA}[c]
    \tilde{\thetav}
    =
    \Proj \tilde{\upsilonv}=\Proj\argmax_{\upsilonv \in \Ups}\LL(\upsilonv)=\argmax_{\thetav \in \Theta} 
        \max_{\Thetathetav{\thetav}} \LL(\upsilonv) =\argmax_{\thetav \in \Theta} \Lr(\thetav).
\label{ttPitu}
\end{EQA}    
The Gauss-Markov Theorem claims the efficiency of such procedures for linear 
Gaussian models and a linear mapping \( \Proj \), and the famous Fisher result extends 
it in the asymptotic sense to the general situation under some regularity conditions. 
The Wilks phenomenon describes the limiting distribution of the likelihood ratio test 
statistic \( T \) which is also called 
the \emph{semiparametric excess}: 
\( T \eqdef 2\bigl\{ \Lr(\tilde{\thetav}) - \Lr(\thetavs) \bigr\} \).
It appears that the distribution of this test statistic is nearly chi-square  
\( \chi^{2}_{\dimp} \) as the samples size grows, \cite{Wilks}:
\begin{EQA}[c]
    T
    \eqdef
    2\bigl\{ \max_{\upsilonv \in \Ups} \LL(\upsilonv) 
    - \max_{\Thetathetav{\thetavs}} \LL(\upsilonv)\bigr\}
    =
    2\bigl\{ \Lr(\tilde{\thetav}) - \Lr(\thetavs) \bigr\}
    \tow 
    \chi^{2}_{\dimp}.
\label{Twilksd}
\end{EQA}    
In particular, the limit distribution does not depend on the particular model 
structure and on the full dimension of the parameter \( \upsilonv \), only the 
dimension of the target matters.
The full parameter dimension can be even infinite under some upper bounds on its 
total entropy. 

The \emph{local asymptotic normality} (LAN) approach by Le Cam leads to the most general 
setup in which the Wilks and Fisher type results can be established. 
However, the classical theory of semiparametric estimation faces serious 
difficulties when the dimension of the nuisance parameter becomes large or infinite. 
The LAN property yields a strong local approximation of the log-likelihood of the full 
model by the log-likelihood of a linear Gaussian model, and this property is only 
validated in a root-n neighborhood of the true point. 
The non- and semiparametric cases require to consider larger neighborhoods where the 
LAN approach is not applicable any more. 
A proper extension of the Wilks and Fisher result to the case of a growing or infinite nuisance 
dimension is quite challenging and involves special constructions like a pilot 
consistent estimator of the target,
a hardest parametric submodel as well as some power tools of the empirical process 
theory; see \cite{MurphyvanderVaart} or \cite{Kosorok} for a comprehensive presentation. 

The recent paper \cite{SP2011} offers a new look at the classical LAN theory. 
The key steps are a local quadratic bracketing for the log-likelihood process and some 
concentration results for its stochastic component. 
The results can be stated for finite samples and do not involve any asymptotic consideration.
It is also shown that many corollaries of the LAN property like Fisher and Wilks exansions only rely on these two facts.
The bracketing idea of \cite{SP2011} is to build two different 
quadratic processes such that the original log-likelihood can be sandwiched between 
them up to a small error. 
This paper offers another approach based on the local linear approximation of the gradient 
of the log-likelihood process. 
This allows to improve the error term of the Fisher and Wilks expansion.
In particular, we obtain a surprising result that the error term in the Fisher expansion can be by factor 
\( \sqrt{\dimtotal} \) smaller than the similar error term in the Wilks Theorem.
In the semiparametric problem with a fixed dimension of the target parameter,
both Fisher and Wilks results apply up to an error \( \dimtotal/\nsize^{1/2} \).
This yields the critical parameter dimension \( \dimtotal= o(\nsize^{1/2}) \).
Another advantage of the new method is that a 
version of the Wilks and Fisher Theorem can be obtained in a quite general semiparametric setup
avoiding any special construction like ``the hardest parametric submodel''.
In addition, it yields the optimal quality of the Fisher and Wilks expansions under the imposed conditions which is confirmed by a specific counter-example. 

For the further presentation we briefly outline the basic steps of the analysis.
Introduce for \(\upsilonv \in \Ups\) and \(\upsilonvs\in\Ups\) as defined in \eqref{eq: definition of full target} the log-likelihood ratio process
\begin{EQA}[c]
    \LL(\upsilonv,\upsilonvs)
    =
    \LL(\upsilonv) - \LL(\upsilonvs) .
\end{EQA}
An important step of our approach is a deviation bound for the MLE \(\tilde{\upsilonv}\in\Ups\).
Given some \(\xx>0\), we define the radius \( \rups = \rups(\xx)>0\) by
\begin{EQA}[c]
    \P(\tilde{\upsilonv} \in \Upss(\rups)) 
    \ge 
    1 - \ex^{-\xx},
\label{devboundupss}
\end{EQA}
where \( \Upss(\rr)\) is a ball of radius \(\rr>0\) in the intrinsic semi-metric corresponding to the process \( \LL(\upsilonv) \). 
We give conditions that ensure that the value \( \rups^{2}(\xx) \) grows almost linearly with 
\( \xx \).
See Section~\ref{sec: conditions} for a precise formulation. 
The second key step is to bound for \(\rr>0\) the approximation error
\begin{EQA}
\label{eq: fundamental linear approximation}
	\UP(\upsilonv)
	& \eqdef &
	\DFc^{-1} \bigl\{ 
		\nabla \LL(\upsilonv) - \nabla \LL(\upsilonvs) - \DFc^{2} \, (\upsilonv - \upsilonvs) 
	\bigr\},
\end{EQA}
where \(\DFc^{2}=\nabla^{2}\E \LL(\upsilonvs)\) is the full information matrix in the model. 
Section~\ref{sec: local approx full space} provides the following bound 
on a set of probability at least \( 1 - \CONST \ex^{-\xx} \):
\begin{EQA}[c]
\label{UPspreadse}
	\sup_{\upsilonv\in\Upss(\rr)}\|\UP(\upsilonv)\|
	\le 
	\Excgr(\rr,\xx) ,
\end{EQA}
where \(\Excgr(\rr,\xx)\) is a small error. 
In combination with the deviation bound \eqref{devboundupss} and 
the identity \( \nabla \LL(\tilde{\upsilonv}) = 0 \),
this allows to derive the following Fisher expansion for the full parameter \( \upsilonv \) inequality:
\begin{EQA}[ccc]
\label{DFspreadse}
	\| \DFc ( \tilde{\upsilonv} - \upsilonvs) - \xiv \|
	& \le & 
	\Excgr(\rr,\xx),
\end{EQA}
where \( \xiv \eqdef \DFc^{-1} \nabla \LL(\upsilonvs) \).
Projecting down on the \( \thetav \)-subspace yields Fisher and Wilks type expansions: 
on a set of dominating probability
\begin{EQA}[ccc]
\label{eq: Fisher in intro}
    \bigl\| \DPr \bigl( \tilde{\thetav} - \thetavs \bigr) - \xivr \bigr\|
	& \leq &
	\Excgr(\rups,\xx),
    \\
\label{eq: Wilks in intro}
    \bigl|\Lr(\tilde{\thetav}) - \Lr(\thetavs) - \| \xivr \|^{2}/2\bigr|
    & \leq &  
	\CONST (\dimp+\xx) \, \Excgr(\rups,\xx)	.
\end{EQA}
The precise definitions of the random \( \dimp \)-vector \( \xivr \) and of the 
symmetric \( \dimp \times \dimp \)-matrix \( \DPr^{2} \) is given below in the next section.
In the case of correctly specified i.i.d models \( \DPr^{2} \) is the covariance matrix of the efficient influence function; see \cite{Kosorok}. 
The random vector \( \xivr \) satisfies \( \E \xivr = 0 \) and 
\( \E \| \xivr \|^{2} \asymp \dimp\) and \(\CONST>0\) is a constant independent of \(\xx>0\) and full dimension \(\dimtotal\). 
Moreover, general deviation bounds for the deviation of quadratic forms from \cite{SP2011} applies to 
\( \| \xivr \|^{2} \) (see Section~\ref{ap: deviation of quadratic forms} for details). 
In the case of a correct model specification they resemble the ones of a chi-square random variable with 
\( \dimp \) degrees of freedom, and the result \eqref{eq: Fisher in intro} can be viewed as an extension of the Wilks phenomenon. 
Under general identifiability conditions, the radius \( \rups \) can be fixed by 
\( \rups^{2} = \CONST_{1}(\dimtotal+\xx) \) for a fixed constant \( \CONST_{1} \) to ensure the concentration property \eqref{devboundupss}.
With this choice of \( \rr \), in the important i.i.d. case, 
the error term \( \Excgr(\rups,\xx) \) can be bounded by 
\( \CONST (\dimtotal+\xx)/\sqrt{\nsize} \).
The results \eqref{eq: Fisher in intro} and \eqref{eq: Wilks in intro} are nonasymptotic and hold true even under model misspecification. 

The proposed approach does not assume that the profile is consistent but gives conditions that ensure the right concentration behavior. Simply assuming that the profile is consistent can be 
even misleading in our setup because this would separate local and global considerations. 
This paper attempts to figure out a list of conditions ensuring global concentration and 
local expansion at the same time. 
This particularly allows to address the crucial question of the largest dimension of the nuisance 
parameter for which the Wilks and Fisher expansions still hold. 
For instance, it appears in regular i.i.d. settings that the condition \( {\dimtotal}^{2} \ll \nsize \) is sufficient where \( \dimtotal \) is the full parameter dimension. 
We present an example that shows that these constraints are critical on the class of considered models. 
It is of interest to compare our statements with the existing literature on the growing parameter asymptotics. 
We particularly mention \cite{Mammen1989,Mammen1993,Mammen1996} and a series of papers by S. Portnoy, see e.g. \cite{Portnoy1984,Portnoy1985,Portnoy1986}.
The typical dimensional asymptotic is \( \dimtotal = o(n^{1/2}) \), which corresponds to our results. 
For some particular special problems and examples the condition on parameter dimension can be relaxed to 
\( \dimp = o(n^{3/2}) \); see \cite{Portnoy1985}.
However, the results are mainly limited to linear or generalized linear regression with independent observations and heavily use the model structure. 
To the contrary, our results apply in a rather general situation and deliver some useful information even in the case when the model is misspecified.
We begin by developing the results for the case that 
the full parameter space \( \Ups \) is a subset of the Euclidean space of 
dimension \( \dimtotal \). 
In Section~\ref{sec: infinite dimensions} we will 
exemplify how to extend our approach to the case when \( \upsilonv \) is a 
functional parameter using the so called sieve approach; see e.g. \cite{ShenShi2005}. 
The present paper combines the sieve approximation idea and 
the finite sample Fisher and Wilks results under a possibly misspecified model. 


The paper is organized as follows. Section~\ref{SWilksFisher} introduces the objects and tools of the analysis and 
collects the main results including an extension of the Wilks Theorem, concentration 
properties of the profile estimator and the construction of confidence sets for the 
``true'' parameter \( \thetavs \). 
Section~\ref{Siidsemi} explains how the results translate to the case of i.i.d. samples and how the approach allows to obtain asymptotic efficiency of the profile estimator in this setting.
Section~\ref{sec: example} presents an example that shows that the ratio \({\dimtotal}^{2}/n\to 0\) is critical to obtain the Wilks phenomenon and the Fisher expansion on the class of models that satisfy the conditions of Section~\ref{sec: conditions}.
%
%
Section~\ref{sec: infinite dimensions} discusses how the results can be extended to the case with the infinite full 
dimension via the sieve approach. We present further conditions on the correlation structure of the full gradient \(\nabla\LL(\upsilonvs)\in\ \BanX\) to also treat the bias.
Section~\ref{sec: application to single index model} briefly outlines how the approach can be employed to derive the main results in the context of single index modeling and which ratio of full dimension to sample size is sufficient in that context. The details of this section can be found in \cite{Andresensingleindex}.
The appendix collects the conditions and the proofs of the main results.

\section{Main results}
\label{sec: main results} 
This section presents our main results on the semiparametric profile estimator
%
which include the Wilks expansion of the profile maximum likelihood \(\Lr(\tilde\thetav,\thetavs)\in\R\) and the Fisher expansion 
of the profile MLE \( \tilde{\thetav}\in\R^{\dimp} \).

Most of results are stated in a finite sample setup for just one fixed sample. 
As we are also interested in understanding what happens if 
\emph{the full dimension} \( \dimtotal \) becomes large 
we also consider a specification of the general finite sample results to 
an asymptotic setup with \( \dimtotal = \dimn \),
where \( \nsize \) denotes the asymptotic parameter, e.g.  
the sample size with \( \nsize \to \infty \). 
%
Our results apply even if the target parameter 
\( \thetav \in\R^{\dimp}\) is also of growing dimension. 
The dimension \( \dimp \) can be of order \( \dimtotal \).
Even the case with a full dimensional target and low dimensional nuisance is 
included as well.

\subsection{Conditions}
\label{sec: conditions}
This section collects the conditions imposed on the model.
We start with the case of a parametric model with a finite dimensional parameter. 
Then explain two new conditions that arise in the case when the nuisance is infinite dimensional.

Let the full dimension of the problem be finite, that is, \(\dimtotal < \infty \).
Our conditions involve two symmetric positive definite \( \dimtotal\times\dimtotal \) matrices, 
the information matrix \( \DFc^{2} \) and the covariance \( \VFc^{2} \),
and a central point \(\upsilonvd\in\R^{\dimtotal}\).
In typical situations for \( \dimtotal < \infty \), 
one can set \(\upsilonvd=\upsilonvs\) where \( \upsilonvs \) is the ``true point'' from 
\eqref{eq: definition of full target}. 
The matrices \( \DFc^{2} \) and \( \VFc^{2} \) can be defined as follows:
\begin{EQA}
	\VFc^{2}
	&=& 
	\Cov\bigl(\score \LL(\upsilonvd)\bigr), 
	\quad
	\DFc^{2}
    =
    - \nabla^{2} \E \LL(\upsilonvd) .
\end{EQA}
Here and in what follows we implicitly assume that the log-likelihood function
\( \LL(\upsilonv) \colon \R^{\dimtotal}\to \R \) is sufficiently smooth in \( \upsilonv \in \R^{\dimtotal}\), 
\( \nabla \LL(\upsilonv)  \in \R^{\dimtotal}\) stands for the gradient and 
\( \nabla^{2} \E \LL(\upsilonv)  \in\R^{\dimtotal\times\dimtotal}\) for the Hessian of the expectation 
\( \E \LL: \R^{\dimtotal}\to \R  \) at \( \upsilonv \in\R^{\dimtotal}\).
It is worth mentioning that such defined matrices \( \DFc^{2} \) and 
\( \VFc^{2}\) are equal if the model \( \Yv \sim \P_{\upsilonvs} \in (\P_{\upsilonv}) \)
is correctly specified and sufficiently regular; see e.g. \cite{IH1981}.


In the context of semiparametric estimation, it is convenient to represent
the information and the covariance matrices in the block form:
\begin{EQA}
    \DFc^{2}
    =
    \left( 
      \begin{array}{cc}
        \DPc^{2} & \Ac \\
        \Ac^{\T} & \HHc^{2} \\
      \end{array}  
    \right),
    & \quad &
    \VFc^{2}
    =
    \left( 
      \begin{array}{cc}
        \VPc^{2} & \Bc
        \\
        \Bc^{\T}  & \VHc^{2} 
        \\
      \end{array}  
    \right) .
\label{DFcseg0}
\end{EQA}
First we state an \emph{identifiability condition}.

\begin{description}
  \item[\( (\bb{\AssId}) \)] 
      There is a constant \( \fis > 0 \) such that 
\begin{EQA}[c]
    \fis^{2} \DPc^{2} 
    \ge  
    \VPc^{2},
    \quad 
    \fis^{2} \HHc^{2} 
    \ge  
    \VHc^{2} ,
    \quad 
    \fis^{2} \DFc^{2}
    \ge 
    \VFc^{2} .
\label{regularity1}
\end{EQA}   
and it holds for some \( \corrDF < 1 \)
\begin{EQA}[c]
    \| \HHc^{-1} \Ac^{\T} \DPc^{-1} \|_{\infty}
    \leq 
    \corrDF .
\label{regularity2}
\end{EQA}
\end{description}

\begin{remark}
The condition \( (\AssId) \) allows to define an important \( \dimp \times \dimp \) efficient information matrix 
\( \DPrc^{2} \) which is defined by inverting the \( \thetav \)-block the full dimensional matrix 
\( \DFc^{2} \).
The exact formula is given by 
\begin{EQA}
	\DPrc^{2}
	& \eqdef &
	\DPc^{2} - \Ac \HH^{-2} \Ac^{\T}
\label{DPrcdefs}
\end{EQA}
and \eqref{regularity2} ensures that the matrix \( \DPrc^{2} \) is well posed. 
Note that the bounds in \eqref{regularity1} are given with the same constant \( \fis \) 
only for simplifying the notation. One can show that the last bound on \( \DFc^{2} \) follows from 
the first two and \eqref{regularity2} with another constant \( \fis' \)
depending on \( \fis \) and \( \corrDF \) only.
\end{remark}

Using the matrix \(\DFc^{2} \) and the central point \(\upsilonvd\in\R^{\dimtotal}\), 
we define the local set \(\Upss(\rr) \subset \Ups \subseteq \R^{\dimtotal}\) with some
\(\rr\ge 0\):
\begin{EQA}[c]
\label{Upssrrdef}
	\Upss(\rr)
	\eqdef 
	\bigl\{ \upsilonv=(\thetav,\etav)\in\Ups \colon \|\DFc(\upsilonv-\upsilonvd)\|\le \rr \bigr\}.
\end{EQA}
Below we distinguish between local and global conditions.
It is assumed that a value \( \rups \) is fixed which separates the local and global zones. 
This value will be specified below.
The local conditions only describe the properties of the process \( \LL(\upsilonv) \) 
for \( \upsilonv \in \Upss(\rr) \) with some fixed value \( \rr \leq \rups \). 
The global conditions have to be fulfilled on the whole \( \Ups \). 
We start with the local conditions. 
The first condition quantifies the smoothness properties of the expected log-likelihood 
\( \E \LL(\upsilonv) \), while the second and the third conditions help to state a similar local
regularity of the stochastic component \( \zeta(\upsilonv) = \LL(\upsilonv) - \E \LL(\upsilonv) \).

\begin{description}
    \item[\( \bb{(\LL_{0})} \)]
    \textit{
    For each \( \rr \le \rups \), 
    there is a constant \( \rddelta(\rr) \) such that
    it holds on the set \( \Upss(\rr) \):
    }
\begin{EQA}[c]
\label{LmgfquadELGP}
    \bigl\|
       \DFc^{-1} \bigl\{ \nabla^{2}\E\LL(\upsilonv) \bigr\} \DFc^{-1} - \Id_{\dimtotal} 
    \bigr\|
    \le
    \rddelta(\rr).
\end{EQA}

\end{description}

\begin{remark}
This condition describes the local smoothness properties of function \( \E \LL(\upsilonv) \).
In particular, it allows to bound the error of local linear approximation of the gradient 
\( \nabla \E \LL(\upsilonv)\): 
under condition \( (\LL_{0}) \) with
\( \DFc^{2} = - \nabla^{2} \E \LL(\upsilonvs) \), it follows from the second order Taylor
expansion for any \( \upsilonv, \upsilonvc \in \Upss(\rr) \)
\begin{EQA}
    \bigl| \DFc^{-1} \bigl\{ 
		\nabla \E\LL(\upsilonv) - \nabla \E\LL(\upsilonvc) 
	\bigr\} - \DFc \, (\upsilonv - \upsilonvc) \bigr|
	& \leq &
	\rddelta(\rr) \rr . 
\label{EdeltGP}
\end{EQA}
In the proofs we actually only need the condition \eqref{EdeltGP} 
which in some cases  can be weaker than \( (\LL_{0}) \).
\end{remark}

The next two conditions concern with the regularity of the stochastic component 
\( \zeta(\upsilonv) \eqdef \LL(\upsilonv) - \E \LL(\upsilonv) \).
Similarly to \cite{SP2011}, we implicitly assume that the stochastic component \( \zeta(\upsilonv) \)
is a separable stochastic process.

\begin{description}
  \item[\( \bb{(\CS \DFc)} \)]
    There exist constants \( \nunu>0 \) and \( \gm > 0 \) such that for all 
    \( |\mubc| \le \gm \)
\begin{EQA}[c]
    \sup_{\gammav \in \R^{\dimtotal}} \log\E \exp\left\{ 
        \mubc \frac{\langle \nabla \zeta(\upsilonvd),\gammav \rangle}
                   {\| \VFc \gammav \|}
    \right\}
    \le 
    \frac{\nu_{0} ^{2} \mubc^{2}}{2}.
\end{EQA}
\end{description}
\begin{remark}
 The matrix \( \VFc^{2} \) describes the variability of the process \( \LL(\upsilonv) \) around 
the central point \( \upsilonvd \) and in many situations can be set as
\begin{EQA}
\label{VFc2se}
    \VFc^{2}
    & \eqdef &
    \Var \bigl\{ \nabla \LL(\upsilonvd) \bigr\}.
\end{EQA}
\end{remark}

\begin{description}
  \item[\( \bb{(\CS \DF_{1})} \)]
    There exists a constant \( \rhor \le 1/2 \), such that for all \( |\mubc| \le \gm \) 
    and all \( 0 < \rr < \rups \)
\begin{EQA}[c]
    \sup_{\upsilonv,\upsilonvc\in\Upss(\rr)}
    \sup_{\|\gammav\|=1} 
    \log \E \exp\left\{ 
         \frac{\mubc \, \gammav^{\T} \DFc^{-1} 
         		\bigl\{ \nabla\zetav(\upsilonv)-\nabla\zetav(\upsilonvc) \bigr\}}
         	  {\rhor \, \|\DFc (\upsilonv-\upsilonvc)\|}\right\}
    \le 
    \frac{\nuno^{2} \mubc^{2}}{2}.
\end{EQA}

\end{description}

%

%

The global conditions are:

\begin{description}
  \item[\( \bb{(\cc{L}{\rr})} \)] 
     For any \( \rr > \rups\) there exists a value \( \gmi(\rr) > 0 \), 
     such that
\begin{EQA}[c]
    \frac{-\E \LL(\upsilonv,\upsilonvd)}{\|\DFc(\upsilonv-\upsilonvd)\|^{2}}
    \ge 
    \gmi(\rr),
    \qquad
    \upsilonv \in \Upss(\rr).
\end{EQA}

  \item[\( \bb{(\CS\rr)} \)] 
    For any \( \rr \ge \rups \) there exists a constant \( \gm(\rr) > 0 \) such that 
\begin{EQA}[c]
    \sup_{\upsilonv \in \Upss(\rr)} \, 
    \sup_{\mubc \le \gm(\rr)} \, 
    \sup_{\gammav \in \R^{\dimtotal}}
    \log\E \exp\left\{ 
        \mubc \frac{\langle \nabla \zeta(\upsilonv),\gammav \rangle}
        {\|\DFc\gammav\|}
    \right\}
    \le \frac{\nunu^{2} \mubc^{2}}{2}.
\end{EQA}
\end{description}

\begin{remark}
The global conditions are the basis for Theorem~\ref{TLDse} and ensure the large deviation property 
\eqref{devboundupss} of \(\tilde{\upsilonv} \) for a properly selected \(\rups\).
In many particular situations, these general conditions can be relaxed.
For instance, in the model with i.i.d. observations, Theorem 5.3 of \cite{IH1981} might serve as a tool. 
The required conditions can be substantially weakened to upper and lower bounds on the Hellinger distance between models for distinct parameters. 
We follow the general way of \cite{SP2011} because it allows to address
possible model misspecification and finite samples.
Note, however, that \( (\cc{L}{\rr}) \) and \( (\CS\rr) \) can be substituted with any other
conditions that describe the value \(\rups\) ensuring \eqref{devboundupss}.
\end{remark}

%

\begin{remark}
We  briefly comment how restrictive the imposed conditions are.
Our conditions on the regularity and smoothness of the log-likelihood process \( \LL(\upsilonv) \)
in terms of the second or even third derivative are stronger than usually required;
cf. Chapters 1,2 in \cite{IH1981}.
Note however, that we do not require that \( \LL(\upsilonv) \) is the true log-likelihood.
It comes from a parametric family chosen by a statistician.
For typical examples, such a family possesses required regularity. 
In particular,
\cite{SP2011}, Section 5.1, considered in details the i.i.d. case and presented some mild 
sufficient conditions on the parametric family which imply the above general conditions.

Conditions \( (\CS \DF_{0}) \), \( (\CS \DF_{1}) \) and \( (\CS \rr) \) requires 
some exponential moments of the observations (errors). 
Usually one only assumes some finite moments of the errors; cf. \cite{IH1981}, Chapter~2.
Our condition is a bit more restrictive but it allows to obtain some finite sample 
bounds. 
Conditions \( (\CS\rr) \) with \( \gm(\rr) \equiv \gm > 0 \) and \( (\cc{L}\rr) \) 
with \( \gmi(\rr) \equiv \gmi > 0 \)
are easy to verify if the parameter set 
\( \Ups \) is compact and the sample size \( \nsize \) exceeds \( \CONST \dimp \) 
for a fixed constant \( \CONST \).
It suffices to check a usual identifiability condition that the value 
\( \E \LL(\upsilonv, \upsilonvs) \) does not vanish for \( \upsilonv \neq \upsilonvs \).

The regression and generalized regression models are included as well; cf. 
\cite{Gh1999,Gh2000} or \cite{Kim2006}.
\cite{SP2011}, Section 5.2, argued that \( (E\!D_{2}) \) is automatically fulfilled for 
a generalized linear model, while \( (E\!D_{0}) \) requires that 
regression errors have to fulfill some exponential moments conditions. 
If this condition is too restrictive and a more stable (robust) estimation procedure is 
desirable, one can apply the LAD-type contrast leading to median regression. 
\cite{SP2011}, Section 5.3, showed for the case of linear median regression that 
all the required conditions are fulfilled automatically if the sample size \( \nsize \) 
exceeds \( \CONST \dimtotal \) for a fixed constant \( \CONST \).
\cite{SpWe2012} applied this approach for local polynomial quantile regression.
\cite{zaitsev2013properties} applied the approach to the problem of regression with Gaussian 
process where the unknown parameters enter in the likelihood in a rather complicated 
way. 
%
\end{remark}

\subsection{Wilks and Fisher expansions}
\label{SWilksFisher} 

This section states the main results in a finite dimensional framework in form of
probabilistic upper bounds on the error of estimation.
The results involve the quantile variable \( \xx \).
If this quantity enters in the bound, it means that this bound is fulfilled with probability at least
\( 1 - \CONST \ex^{-\xx} \) for a fixed explicit constant \( \CONST \).


First we introduce the main elements of the approach. 
Let the \emph{information matrix} \( \DFc^{2}\in\R^{\dimtotal\times\dimtotal}\) be from \( (\LL_{0}) \), 
and the \emph{score covariance matrix} \( \VFc^{2}\in\R^{\dimtotal\times\dimtotal}\) be from 
\( (\CS \DFc) \).
Introduce the \emph{misspecification matrix} \( \BB\in\R^{\dimtotal\times\dimtotal}\)
given by the famous sandwich formula; see \cite{huber1967}:
\begin{EQA}
\label{DFc2se}
    \BB
    &=&
    \DFc^{-1} \VFc^{2} \DFc^{-1}.
\end{EQA}
%
In the case of correct model specification with \( \DFc^{2} = \VFc^{2} \), 
the \emph{sandwich matrix} \( \BB \) becomes the identity:
\( \BB = \Id_{\dimtotal} \).

For the semiparametric \( (\thetav,\etav) \)-setup, we
consider the block representation of the vector \( \score \eqdef \nabla \LL(\upsilonvs) \) 
and of the matrices \( \DFc^{2}, \VFc^{2}\) from \eqref{DFc2se}: 
\begin{EQA}
    \score
    =
    \left( 
      \begin{array}{c}
        \score_{\thetav} \\
        \score_{\etav} 
      \end{array}  
    \right),
    & \quad
    \DFc^{2}
    =
    \left( 
      \begin{array}{cc}
        \DPc^{2} & \Ac \\
        \Ac^{\T} & \HHc^{2} 
      \end{array}  
    \right),
    & \quad
	\VFc^{2}=\left( 
      \begin{array}{cc}
        \VPc^{2} & \Bc \\
        \Bc^{\T} & \VH^{2} 
      \end{array}  
    \right).
\label{DFcseg}
\end{EQA}  
\def\VPr{\breve{\VP}}
\def\BBr{\breve{\BB}}
Define also the \( \dimp \)-vectors \( \scorer_{\thetav} \) and \( \xivr \in \R^{\dimp} \)
\begin{EQA}
	\scorer_{\thetav}
    &=&
    \score_{\thetav} - \Ac \HHc^{-2} \score_{\etav},
    \quad
    \xivr
    \eqdef
    \DPr^{-1} \scorer_{\thetav} \, .
\label{scorerdef}
\end{EQA}
and \( \dimp \times \dimp \) matrices \( \DPr^{2} ,\VPr^{2}, \BBr\) as
\begin{EQA}
    \DPr^{2}
    &=&
    \DP^{2} - \Ac \HHc^{-2} \Ac^{\T}, 
    \\
    \VPr^{2}
    &=&
    \Cov(\scorer_{\thetav})
    =
    \VPc^{2} - 2 \Bc \HHc^{-2} \Bc^{\T} + \Bc \HHc^{-2} \VH^{2} \HHc^{-2} \Bc^{\T}, \quad
    \\
    \BBr
    &=&
    \DPr^{-1}\VPr^{2}\DPr^{-1}.
\label{DFse02}
\end{EQA}
The random variable \(\scorer_{\thetav}\in\R^{\dimp} \) is related to the efficient influence function in semiparametric estimation and the matrix \(\DPr^{2}\in \R^{ \dimp \times \dimp}\) equals its covariance in the case of correct specification. 
Finally we introduce \( \tilde{\upsilonv}_{\thetavs}\in\Ups \), which maximizes \( \LL(\upsilonv,\upsilonvs) \) subject to
\( \Proj \upsilonv = \thetavs \):
\begin{EQA}[c]
\label{tuthLLuus}
    \tilde{\upsilonv}_{\thetavs} \eqdef (\thetavs,\tilde \etav_{\thetavs})
    \eqdef 
    \argmax_{\Thetathetav{\thetavs}} 
    \LL(\upsilonv,\upsilonvs).
\end{EQA}
The Wilks expansion below claims that the profile maximum likelihood
\( \Lr(\tilde{\thetav},\thetavs) \eqdef \Lr(\tilde{\thetav}) - \Lr(\thetavs) \)
 can be approximated by a quadratic form \( \| \xivr \|^{2}/2 \) with 
 \( \xivr = \DPr^{-1} \scorer_{\thetav} \).
In the case of correct model specification  
the deviation properties of the quadratic form 
\( \| \xivr \|^{2} = \| \DPr^{-1} \scorer_{\thetav} \|^{2} \) are essentially 
the same as those of a chi-square random variable with \( \dimp\in\N \) degrees of freedom;
see Proposition~\ref{theo: dev bounds quad forms} in the appendix.
Moreover, in the correctly specified i.i.d setup the vector \( \xivr \) is asymptotically standard 
normal; see Section~\ref{Siidsemi}.
In the general case, the behavior of the quadratic 
form \( \| \xivr \|^{2} \) depends on the characteristics of the matrix
\( \BBr \eqdef \DPr^{-1} \VPr^{2} \DPr^{-1} \).
More precisely,
the presented finite sample results involve the upper quantile function \( \zz(\xx,\BBr) \) 
of this quadratic form ensuring 
\begin{EQA}
	\P\bigl( \| \xivr \| > \zz(\xx,\BBr) \bigr) 
	& \leq &
	2 \ex^{-\xx} ; 
\label{PxivrzzBBr}
\end{EQA}
see again Proposition~\ref{theo: dev bounds quad forms}.
One can use the bound \( \zz^{2}(\xx,\BBr) \asymp \dimp + \xx \) in most of situations.

Further important objects in our results are the central point \( \upsilond \) and the local radius \( \rups \) that ensures under conditions 
\({(\cc{L}{\rr})} \) and \({(\CS\rr)} \) that the MLE \( \tilde{\upsilonv} \) 
lies with a high probability in the local vicinity 
\( \Upss(\rups) = \bigl\{ \upsilonv\in\Ups \colon \|\DFc(\upsilonv-\upsilonvd)\|\le \rups \bigr\} \);
see \eqref{devboundupss}.
Theorem~\ref{TLDse} in Section~\ref{sec: large deviations} states such a result for the value 
\(\rups = \rups(\xx) \approx \CONST\sqrt{\xx+\dimtotal}>0\).

Define the \emph{semiparametric spread} \(\Excgr(\rr,\xx)>0\) as
\begin{EQA}[c]
    \Excgr(\rr,\xx)
    \eqdef
    \left\{ \rddelta(\rr) + 6 \, \nuno \, \zzQ(\xx,\dimtotal) \, \rhor \right\} \rr,
\label{eq: def of diamond rr}
\end{EQA}    
where \( \rddelta(\rr) \) is shown in the condition \({ (\LL_{0})} \), 
the constants \(  \rhor \), \(\nuno \) are from condition \({ (\CS \DF_{1})} \) 
in Section~\ref{sec: conditions}. 
The value \( \zzQ(\xx,\dimtotal) \) is related to the entropy of the unit ball
in a \( \R^{\dimtotal} \)-dimensional Euclidean space and 
one can apply \(\zzQ(\xx,\dimtotal)\cong \sqrt{\xx+4\dimtotal}\)
 for moderate choice of \(\xx>0\) (see Appendix~\ref{ap: bound for norm of gradient}). 
The value \( \Excgr(\rr,\xx) \) measures the quality of a linear approximation to \(\nabla\LL(\upsilonv)-\nabla\LL(\upsilonvs)\) in the local vicinity \( \Upss(\rr) \); see \eqref{DFspreadse}. 
Our results become accurate if \(\Excgr(\rups,\xx)\) is small, where \(\rups>0\) is chosen to ensure 
the large deviation result \eqref{devboundupss}.
The spread will be evaluated in the i.i.d. case in Section~\ref{Siidsemi} below.

%
\begin{theorem}
\label{theo: main theo finite dim}
Assume \({(\CS \DFc)} \), \({(\CS \DF_{1})} \), \({(\LL_{0})}\), and \({(\AssId)}\) with 
a central point \(\upsilonvd \) 
and some matrices \( \DFc^{2} \) and 
\( \VFc^{2} \). 
Fix a radius \(\rups>0\) such that 
\(\P(\tilde{\upsilonv},\tilde{\upsilonv}_{\thetavs}\in\Upss(\rups))\ge 1-\ex^{-\xx}\). 
Then it holds on a set \(\Omega(\xx)\subseteq\Omega \) of probability at least
\( 1-6\ex^{-\xx} \) for the profile MLE \( \tilde{\thetav} \) from \eqref{Lrthetav} and \eqref{ttPitu}
\begin{EQA}
	\bigl\| 
        \DPr \bigl( \tilde{\thetav} - \thetavs \bigr) 
        - \xivr 
    \bigr\|
    &\le& 
    \Excgr(\rups,\xx) ,
\label{eq: Fisher in main theo}
	\\
    \bigl| 2 \Lr(\tilde{\thetav},\thetavs) - \| \xivr \|^{2} \bigr|
    &\le&
    2\rr_{\dimp} \Excgr(\rups,\xx),
\label{eq: Wilks in main theo}
\end{EQA}
where the spread \( \Excgr(\rups,\xx) \) is defined in \eqref{eq: def of diamond rr} and
\begin{EQA}[c]
 	\rr_{\dimp}
	\eqdef
	\left(3+\frac{1+\corrDF}{1-\corrDF}\right) \Excgr(\rups,\xx)
	+ \frac{1+\corrDF}{1-\corrDF} \zz(\xx,\BBr).
\end{EQA}
\end{theorem}

\begin{remark}
{
In the classical finite dimensional case, a usual choice for the central point \( \upsilonvd \) is
\(\upsilonvd=\upsilonvs = \argmax_{\upsilonv \in \Ups} \E \LL(\upsilonv)\) 
and one can define the matrices \( \DFc^{2} \) and \( \VFc^{2} \) as \( \DFc^{2} = - \nabla^{2} \E \LL(\upsilonvs) \) and 
\( \VFc^{2} = \Cov\bigl(\score \LL(\upsilonvs)\bigr) \). 
However, for the sieve semiparametric problem below, we use another definition related to the infinite 
dimensional model.} 
\end{remark}

\begin{remark}
\label{Rsemib3}
The profile maximum likelihood process \( \Lr(\thetav) \) can be used for defining 
the likelihood-based confidence sets of the form
\begin{EQA}[c]
    \CS(\zz) 
    = 
    \{ \thetav:  \Lr(\tilde{\thetav},\thetav) \le \zz \}
\label{Cszzse}
\end{EQA}    
for some \( \zz > 0 \). 
The bound \eqref{eq: Wilks in main theo} helps to evaluate the coverage probability 
\( \P\bigl( \thetavs \notin \CS(\zz) \bigr) \) in terms of deviation probability for 
the quadratic form \( \| \xivr \|^{2} \); cf. Corollary~3.2 in \cite{SP2011}.
\end{remark}

\begin{remark}
One can use the expansion \eqref{eq: Fisher in main theo} for describing the concentration 
probability for elliptic sets 
\begin{EQA}[c]
    \CA(z) 
    =
    \bigl\{ 
        \thetav: \| \DPr (\thetav - \thetavs) \|
        \le 
        z
    \bigr\} ;
\label{Cazzloc}
\end{EQA}
cf. Corollary 3.5 in \cite{SP2011}.
\end{remark}

\subsection{Large deviation bounds}
\label{sec: large deviations}
In this section we want to present a way to determine a value \(\rups>0\) such that the full MLE 
\( \tilde{\upsilonv} \in\R^{\dimtotal}\) belongs to the local vicinity \(\Upss(\rups) \subset \R^{\dimtotal}\) with dominating probability. 
It is important to note that Theorem~\ref{TLDse} is one particular approach 
which could be replaced by any other proper technique.

As a first step we adopt the upper function approach from \cite{SP2011}; cf. 
Corollary 4.4 therein. 
Again the constants \( \gm(\rr) \) and \( \gmi(\rr) \) are introduced in 
Section~\ref{sec: conditions}.

\begin{theorem}[\cite{SP2011}, Theorem~4.1]
\label{TLDse}
Suppose \( (\CS\rr) \) and \( (\cc{L}\rr) \) with \( \gmi(\rr) \equiv \gmi \).
If for a fixed \( \rups \) and any \( \rr \ge \rups \), the following conditions are fulfilled:
\begin{EQA}
    1 + \sqrt{\xx + 2\dimtotal} 
    & \le &
    3 \nu_{\rr}^{2} \gm(\rr)/\gmi ,
    \\
    6 \nu_{\rr} \sqrt{\xx + 2\dimtotal}
    & \le &
    \rr\gmi ,
\label{cgmi2rrc}
\end{EQA}
then 
\begin{EQA}
	\P(\tilde{\upsilonv},\tilde{\upsilonv}_{\thetavs}\in\Upss(\rups))
	& \ge &
	1-\ex^{-\xx}.
\label{PLDttuttut}
\end{EQA}
\end{theorem}

\begin{remark}
\label{remark: bound for rups}
The condition \eqref{cgmi2rrc} helps to understand which \( \rups >0\) ensures 
prescribed concentration properties of \( \tilde{\upsilonv}\in \R^{\dimtotal} \) and 
\( \tilde{\upsilonv}_{\thetavs} \in \R^{\dimtotal}\).
Namely, if \( \gm(\rr)>0 \) is large enough, then \eqref{cgmi2rrc} follows from the 
bound 
\begin{EQA}[c]
    \rr_0
    \ge 
    6 \gmi^{-1} \nu_{\rr} \sqrt{\xx + \dimtotal} .
\label{rupslbse}
\end{EQA}    
\end{remark}

The upper function approach of showing the consistency for an M-estimator can be rather rough
and the bound \eqref{rupslbse} could lead to rather large values of \(\rups>0\). 
As the obtained value \(\rups>0\) enters into the error term \(\Excgr(\rups,\xx) >0\) of 
Theorem~\ref{theo: main theo finite dim} it is desirable to obtain a general refined bound for 
\( \rupf\le \rups\) that still ensures that 
\(\P(\tilde{\upsilonv}\in\Upss(\rupf))\ge 1-\CONST\ex^{-\xx}\) with a small constant \(\CONST>0\). 
Such an improvement is possible if the spread \( \Excgr(\rups,\xx) \) describing the quality of local 
approximation is not too big, namely, if \(\Excgr(\rups,\xx)\le \zz(\xx,\BB)/(1-\corrDF) \).

\begin{proposition}
\label{lem: large dev refinement}
Assume the conditions of Theorem~\ref{theo: main theo finite dim}. 
Let \(\rups>0\) be such that \eqref{PLDttuttut} holds and define the radius \( \rupf \le \rups \)
\begin{EQA}[c]
\label{eq: def or refined radius}
	\rupf 
	\eqdef 
	\zz(\xx,\BB)/(1-\corrDF)+\Excgr(\rups,\xx)\wedge \rups.
\end{EQA}
Then the result of Theorem~\ref{theo: main theo finite dim} continues to apply with the error term 
\(\Excgr(\rupf,\xx)\) in place of \(\Excgr(\rups,\xx)>0\).
\end{proposition}
%
%

\subsection{The i.i.d. case}
\label{Siidsemi}
Here we briefly discuss the implications of our general results to the case 
with \( \Yv = (Y_{1},\ldots,Y_{\nsize})^{\T} \), 
where observations \( Y_{i} \) are i.i.d. from a measure \( P \) and \( \nsize \) is the sample size.
The parametric assumption means \( P = P_{\upsilonvs} \in (P_{\upsilonv}, \upsilonv \in \Upsilon\subseteq \R^{\dimtotal}) \) 
for a given parametric family \( (P_{\upsilonv}) \).
We assume that \( (P_{\upsilonv}) \) obeys the regularity conditions 
listed in Section~5.1 of \cite{SP2011}.
By \( \ell(y,\upsilonv) \in\R\) we denote the log-density of \( P_{\upsilonv} \) w.r.t. 
some dominating measure \( \Pdom \).
For simplicity of comparison with the classical results 
we do not discuss the model misspecification issue, i.e. 
assume that the parametric assumption is correct.
However, an extension to the case of a misspecified model is straightforward.
Assume that the likelihood for a single observation \( \ell(y,\upsilonv)\in\R \) satisfies the conditions in Section~\ref{sec: conditions} with some \(\nunu\), \(\rhorb \), \(\rddelta(\rr)= \rddeltab\rr \), 
\(\gmi(\rr)=\gmi^{*}\) and \( \gm=\gmone \) for some positive constants \(\nunu, \rhorb ,\rddeltab, \gmi^{*}, \gmone>0\). 
Under these conditions, one can easily check the conditions in Section~\ref{sec: conditions} for 
the full log-likelihood 
\( \LL(\upsilonv) = \sum_{i=1}^n\ell(y_i,\upsilonv)\) with  
\(\rhor= \rhorb \nsize^{-1/2} \), \(\rddelta(\rr)= \rddeltab\rr\nsize^{-1/2} \), 
\(\gmi(\rr)=\gmi^{*}\), and \( \gm=\gmone\nsize^{1/2} \), 
and with \(\DFc^{2} = \VFc^{2} = \IF \) being the Fisher information matrix of the family 
\( (P_{\upsilonv}) \) at the point \( \upsilonvs \); cf. Lemma 5.1 in \cite{SP2011}.
Theorem~\ref{TLDse} yields that
\begin{EQA}
	\P\left(\tilde{\upsilonv}\in\Upss(\rups)\right)\ge 1-\ex^{-\xx}, 
	&\text{ with } &
	\rups(\xx)=6\frac{\nunu}{\gmi^{*}} \sqrt{2\dimtotal+\xx} .
\end{EQA}
Also one obtains
\begin{EQA}[c]
	\Excgr(\rr,\xx)
	= 
	\rr(\rddelta^{*}\rr+6\nunu\rhor^{*}\sqrt{4\dimtotal+\xx})/\sqrt{\nsize}.
\end{EQA}
%
%
To apply Theorem~\ref{theo: main theo finite dim} we further need a version of the identifiability condition 
\({(\AssId)} \) from Appendix~\ref{sec: conditions} on the marginal distribution.
Consider the block representation of the Fisher information matrix \( \IF \):
\begin{EQA}[c]
    \IF
    =
    \left( 
      \begin{array}{cc}
        \IF_{\thetav\thetav} & \IF_{\thetav\etav} \\
        \IF_{\thetav\etav}^{\T} & \IF_{\etav\etav}
      \end{array}  
    \right) .
\label{iiciid}
\end{EQA}    
The required identifiability condition reads as follows:
\begin{description}
  \item[\( (\bb{\assId}) \)] 
      There is a constant \( \corrDF < 1 \) such that 
\begin{EQA}[c]
    \| \iic_{\etav\etav}^{-1/2}\iic_{\thetav\etav}^{\T} \iic_{\thetav\thetav}^{-1/2} \|_{\infty}
    \le 
    \corrDF .
\label{regularity2iid}
\end{EQA}
\end{description}

Also define
\begin{EQA}[c]
    \ifc
    \eqdef
    \IF_{\thetav\thetav} - \IF_{\thetav\etav} \IF_{\etav\etav}^{-1} \IF_{\thetav\etav}^{\T} .
\label{ifcdpcahc}
\end{EQA}
Now Theorem~\ref{theo: main theo finite dim} applies and obtain the following result.

\begin{corollary}
\label{Txiviidse}
Let \( Y_{1},\ldots,Y_{\nsize} \) be i.i.d. \( \P_{\upsilonvs} \) and 
let the likelihood for a single observation \( \ell(y,\upsilonv)\in\R \) satisfy the conditions in Section~\ref{sec: conditions} with \(\nunu=\nunu\), \(\rhor= \rhorb \), \(\rddelta(\rr)= \rddeltab\rr \), \(\gmi(\rr)=\gmi^{*}\) and \( \gm=\gmone \) for some positive constants \(\nunu, \rhorb ,\rddeltab, \gmi^{*}, \gmone>0\).
In addition, assume \( (\bb{\assId}) \); see \eqref{regularity2iid}.
Then we get the Fisher and Wilks results of Theorem~\ref{theo: main theo finite dim} with \(\BB=\Id_{\dimtotal}\), \(\BBr=\Id_{\dimp}\), and
\begin{EQA}[c]
	\Excgr(\rups,\xx)
	\le 
	36\frac{\nunu}{\gmi^{*}}\left(\rddelta^{*}\frac{\nunu}{\gmi^{*}}(\xx+2\dimtotal)^{2} 
	+ \nunu\rhor^{*}\zzQ(\xx,\dimtotal)\sqrt{\xx+2\dimtotal}\right)/ \sqrt{\nsize}.
\end{EQA}
\end{corollary}    

\begin{remark}
The definition of \(\zzQ(\xx,\dimtotal) \) implies for moderate values of \(\xx>0\) that
\begin{EQA}[c]
	\Excgr(\rups,\xx)
	\le
	\CONST_{\zz}(\xx+\dimtotal)/{\sqrt{\nsize}}
\end{EQA}
with some fixed constant \( \CONST_{\zz} \).
The Fisher result \eqref{eq: Fisher in main theo} is meaningful if 
\( \Excgr(\rups,\xx) \) is small yielding the constraint \( \dimtotal \ll \nsize^{1/2} \).
If the target dimension \(\dimp\) is fixed, the same condition is sufficient for the Wilks expansion 
in \eqref{eq: Wilks in main theo}.
However, if the target dimension \(\dimp\) is of order \(\dimtotal\), the constraint for the Wilks theorem becomes \(\dimtotal=o(n^{1/3})\).
\end{remark}

\subsection{Critical dimension}
\label{sec: example}
This section discusses the issue of \emph{critical parameter dimensions} when the full 
dimension \( \dimtotal \) grows with the sample size \( \nsize \).  
We write \( \dimtotal = \dimn \). 
The results of Proposition~\ref{lem: large dev refinement} are accurate if 
the spread function 
\( \Excgr(\rr,\xx) \) from \eqref{eq: def of diamond rr} fulfills
\(\Excgr(\rups,\xx) \le  \zz(\xx,\BB)\) and \(\Excgr(\rupf,\xx)\) is small, with 
\( \rupf = 2\zz(\xx,\BB)/(1-\corrDF) \).
Usually \(\zz(\xx,\BB)\le \CONST \sqrt{\xx+\dimtotal}\) leading to 
\begin{EQA}[c]
\label{eq: reduce spread statement critical}
	\Excgr(\rupf,\xx)
	\asymp 
	\delta(\rupf) \rupf + \rhor \rupf^{2}
	\quad
	\text{ is small for }
	\rupf^{2} \asymp \dimtotal.
\end{EQA}
The critical size of \(\dimtotal\) then depends on the exact bounds on \(\delta(\cdot),\rhor\). 
If \(\delta(\rr)/\rr \asymp \rhor \asymp 1/\sqrt{\nsize}\) (as in Theorem~\ref{Txiviidse}) the condition \eqref{eq: reduce spread statement critical} reads \(\Excgr(\rupf,\xx) \asymp \dimtotal/\sqrt{\nsize}\). 
In other words, one needs that ``\({\dimtotal}^{2}/n\) is small'' to obtain 
an accurate non asymptotic version of the Wilks phenomenon and the Fisher Theorem. 
Similar conclusions were obtain by Portnoy in series of papers on growing dimension in generalized linear models and for natural exponential families, see e.g. \cite{Portnoy1984,Portnoy1985,Portnoy1986}. 
Our results are non-asymptotic and apply to general statistical models under the conditions 
from Section~\ref{sec: conditions}. 
The following example shows that the constraint ``\({\dimtotal}^{2}/n\) is small'' is critical.

%
Consider \(n\in\N\) i.i.d. observations in the model
\begin{EQA}
\Yv_i&=&f(\upsilonv)+\varepsilonv_i,\\
f(\upsilonv)&=&f(\thetav,\etav)\eqdef\left(\begin{array}{l}
		     \thetav \\
		      \etav_1\\
		      \vdots\\
		      \etav_{\dimn-1}\\
                     \end{array}\right)
  +\left(\begin{array}{l}
		     \|\etav\|^{2} \\
		      0\\
		      \vdots\\
		      0\\
                     \end{array}\right)\in\R^{\dimn},
\end{EQA}
with \(\varepsilonv_i\sim \mathcal N(0,I_{\dimtotal})\) and \(\upsilonv=(\thetav,\etav)\in\R\times \R^{\dimn-1}\). Assume that the parameter of interest is \(\thetav\in\R\) and that the true point satisfies \(\upsilonvs=0\in \R^{\dimtotal}\). 

\begin{proposition}
 \label{lem: critical dim in fisher}
Under \(\dimn/\sqrt{\nsize} \to 0\), the Fisher expansion is accurate. 
If \(\dimn/\sqrt{\nsize}\not\to 0\) the profile MLE in the above model is not root-n consistent. 
For \(\sqrt{\nsize}=o(\dimn)\) the root-n bias tends to infinity almost surely. 
Finally, the Wilks phenomenon occurs if and only if \(\dimn=o(\sqrt{\nsize})\).
\end{proposition}

\subsection{Infinite dimensional nuisance}
\label{sec: infinite dimensions}

This section discusses how the approach can be extended to the infinite dimensional 
case. First the basic idea of projecting the infinite dimensional problem down to a 
finite dimensional one is explained. 
Then we prove under bias constraints that the 
projected \emph{sieve} estimator is nearly normal and efficient. 
Finally the approach is illustrated with the single index model as an example. For simplicity we present the case of a Hilbert space. The ideas can be modified to treat the case when 
the nuisance parameter belongs to a Banach space. 

\subsubsection{Sieve approach}
Consider the 
\( (\thetav,\fv) \)-setup with \( \thetav \in \Theta \subseteq \R^{\dimp} \) and
\( \fv \in \BanX \), where \( \BanX \) is an infinite dimensional 
Hilbert space. 
The target parameter \( \thetavs \) can be defined as 
\begin{EQA}[c]
    \thetavs
    =
    \argmax_{\thetav} \sup_{\fv \in \BanX} \E \LL(\thetav,\fv) .
\label{thetavsBanX}
\end{EQA}
Assume that the infinite dimensional 
Hilbert space \( \BanX \) possesses a countable orthonormal basis \( \{\basX_{1},\basX_{2},\ldots\} \subset \BanX\). The vector 
\( \fv\in \BanX  \) admits a unique decomposition in the form 
\begin{EQA}[c]
    \fv
    =
    \sum_{j=1}^{\infty} \eta_{j} \basX_{j}.
\label{fsexpn}
\end{EQA} 
As \( \BanX \) is a Hilbert space \( \eta_{j} = \bigl\langle \fv, \basX_{j} \bigr\rangle \) is the usual Fourier 
coefficient. In the \emph{sieve} approach one assumes that for any \( \dimh\in\N \) a finite set 
\( \basX_{1},\ldots,\basX_{\dimh} \) of elements in \( \BanX \) is fixed and the 
vector \( \fv \) can be approximated by a finite linear combination 
\( \fv_{\dimh}(\etav) \) of the \( \basX_{j} \)'s:
\begin{EQA}[c]
    \fv_{m}(\etav)
    \eqdef
    \sum_{j=1}^{m} \eta_{j} \basX_{j}.
\label{fsdimhse}
\end{EQA}    
We denote the parameter by \( \upsilonv = (\thetav,\etav)\in l^{2} \).
In the following we will need to quantify the accuracy of approximating 
\( \fv \) by \( \fv_{m} \) as \( m \) grows; see condition \( \bb{(bias)} \) below.

Let \( \LL(\thetav,\fv) \) be the log-likelihood in the original model. Define
\begin{EQA}
\LL(\upsilonv)&\eqdef&\LL\left(\thetav,\sum_{j=1}^{\infty}\eta_j\basX_j\right)\\
\upsilonvs&\eqdef&\argmax_{(\thetav,\etav)\in l^{2}}\E\left[\LL\left(\thetav,\sum_{k=1}^{\infty}\eta_j\basX_j\right)\right],
\end{EQA}
and the \( \dimh \)-dimensional sieve approximation \( \LL_{\dimh}(\upsilonv) \) 
of \( \LL(\upsilonv) \) by
\begin{EQA}
    \LL_{\dimh}(\thetav,\etav)
    & \eqdef &
    \LL(\thetav,\fv_{\dimh}(\etav)) ,
    \qquad 
    (\thetav,\etav) \in \Upsilon_{\dimh}\eqdef\{\upsilonv=(\thetav,\etav)\in \R^{\dimtotal}: (\thetav,\fv_\dimh(\etav))\in\Upsilon \}.
\label{LLdimhte}
\end{EQA}    
The corresponding sieve profile estimator \( \tilde{\thetav}_{\dimh} \) and its target 
\( \thetavs_{\dimh} \) for this parametric \( \dimh \)-submodel are defined in the 
usual way:
\begin{EQA}
\label{ttdhttdhs}
    \tilde{\thetav}_{\dimh} 
    & \eqdef &\Proj\tilde{\upsilonv}_{\dimh} \eqdef \argmax_{\upsilonv \in \Upsilon_{\dimh}}
        \LL_{\dimh}(\thetav,\etav),
    \\
    \thetavs_{\dimh} 
    & \eqdef &\Proj\upsilonvs_{\dimh} \eqdef\argmax_{\upsilonv \in \Upsilon_{\dimh}} 
        \E \LL_{\dimh}(\thetav,\etav).
\label{tdhttdhs}
\end{EQA}
The question we are interested in can be formulated as follows:
is \( \tilde{\thetav}_{\dimh} \) a good (efficient) estimator of \( \thetavs \) from 
\eqref{thetavsBanX} under a proper choice of \( \dimh \)?

\subsubsection{Bias constraints and efficiency}
The parametric results obtained in Section~\ref{sec: main results} claim that 
\( \tilde{\thetav}_{\dimh} \in\R^{\dimp}\) estimates well \( \thetavs_{\dimh} \in\R^{\dimp} \) if the spread \(\Excgr(\rups,\xx)>0\) is small. More precisely we have the following: Define for fixed \(\xx>0 \) the value \( \rups >0\) by 
\(\P\big\{\tilde{\upsilonv},\tilde{\upsilonv}_{\thetavs}\in\Upss(\rr_{0,n})\big\}\ge 1-\ex^{\xx}\). Applying Theorem~\ref{theo: main theo finite dim} to 
\( \tilde{\thetav}_{\dimh} \) from \eqref{ttdhttdhs} we find that with probability greater \(1-6 \ex^{-\xx}\)
\begin{EQA}[c]
\label{eq: convergence of estimator for eta}
    \|\DPrp \bigl( \tilde{\thetav}_{\dimh} - \thetavs_{\dimh} \bigr)
    - \xivr_{\dimh}(\upsilonvs_{\dimh})\|\le \Excgr(\rups,\xx),
\label{DPrctthsh}
\end{EQA} 
where \( \upsilonvs_{\dimh} = (\thetavs_{\dimh},\etavs_{\dimh}) 
= \argmax_{\upsilonv} \E \LL_{\dimh}(\upsilonv) \) and 
\begin{EQA}
\DPr^{2}_{\dimh}\eqdef\left(\Proj\DF_{\dimh}^{-2}\Proj^{\T}\right)^{-1}\in\R^{\dimp\times\dimp},&&
\DF^{2}_{\dimh}\eqdef\nabla_{\dimp+\dimh}^{2}\E[\LL(\upsilonvs_{\dimh})]\in\R^{\dimtotal\times\dimtotal},
\end{EQA}
i.e. the derivatives of \(\E[\LL]\) are only taken with respect to the first \(\dimp+\dimh\in\N\) coordinates of \(\upsilonv\in l^{2}\).
This result involves two kinds of bias, one that concerns the difference \(\thetavs_{\dimh} - \thetavs\) and the other the difference between \(\DPrp\in\R^{\dimp\times\dimp}\) and \(\DPr\in\R^{\dimp\times\dimp}\) where
\begin{EQA}
\DPr^{2}\eqdef \left(\Proj\DF^{-2}\Proj^{\T}\right)^{-1}\in\R^{\dimp\times\dimp},&&
\DF^{2}\eqdef\nabla^{2}\E[\LL(\upsilonvs)]\in L(l^{2},l^{2}),
\end{EQA}
i.e. the derivatives of \(\E[\LL]\) are take with respect to all coordinates of \(\upsilonv\in l^{2}\) and the Hessian is calculated in the "true point" \(\upsilonvs\in l^{2}\). The second bias - i.e. bounds for \(\| I-\DPrp^{-1}\DPr^{2}\DPrp^{-1}\|\) - will be neglected for now, as only the operator \(\DPr^{2}_{\dimh}\in\R^{\dimp\times\dimp}\) is available in practice. We will come back to it, when we derive efficiency for the sieve profile estimator \(\tilde \thetav_{\dimh}\in\R^{\dimp}\).

For the first type of bias we impose the following condition:

\begin{description}
\item[\(\bb{(bias)}\)] There exists a  decreasing function \(\alpha:\N\to \R_+\) such that
\begin{EQA}[c]
\|\DPrp(\thetavs_{\dimh} - \thetavs) \|\le \alpha(\dimh).
\end{EQA}
\end{description}

\begin{remark}
This paper focuses on the result~\ref{theo: main theo finite dim} and thus we do not elaborate on approximation theory. But \cite{AAbias2014} presents conditions on the structure of \(\DF: l^{2}\to l^{2}\) and on the sequence \(\etavs\in l^{2}\) that yield \(\bb{(bias)}\) in an adequate way. 
\end{remark}

\medskip
For \(\dimh\in\N\) and \(\rr\ge 0\) define the local set \({\Upsilon}_{0,\dimh}(\rr)\)
\begin{EQA}[c]
{\Upsilon}_{0,\dimh}(\rr)\eqdef \{\upsilonv=(\thetav,\etav)\in\R^{\dimtotal}:\, (\thetav,\fv_{\dimh}(\etav))\in\Upsilon,\, \|\DF_{\dimh}(\upsilonv-\upsilonvs_{\dimh})\|\le \rr\},
\end{EQA}
and \( \tilde{\upsilonv}_{\thetavs,\dimh}\in\Ups \), which maximizes \( \LL_{\dimh}(\upsilonv,\upsilonvs) \) subject to
\( \Proj \upsilonv = \thetavs \):
\begin{EQA}
\label{tuthLLuus}
    \tilde{\upsilonv}_{\thetavs_{\dimh},\dimh}
    & \eqdef &
    \argmax_{\substack{\upsilonv\in\Upsilon \\ \Proj\upsilonv=\thetavs_{\dimh}}} 
    \LL_{\dimh}(\upsilonv,\upsilonvs).
\end{EQA}
Further we represent
\begin{EQA}[c]
\DF_{\dimh}=\left(\begin{array}{cc} \DP^{2} & \A_{\dimh}^{\T}\\
			\A_{\dimh} & \HH^{2}_{\dimh} \end{array}\right),
\end{EQA}

With Theorem~\ref{theo: main theo finite dim} and \(\bb{(bias)}\) we directly get the following corollary:

\begin{corollary}
\label{cor: semi sieve bias}
Assume \(\bb{(bias)}\) and that the conditions \(\bb{ (\CS \DFc)} \), 
\(  \bb{(\CS \DF_{1})} \) and \( \bb{ (\LL_{0})} \) 
from Section~\ref{sec: conditions} are satisfied for all \(\dimh\ge\dimh_0\) for some \(\dimh_0\in\N\) and with \(\DF^{2}=\nabla_{\dimp+\dimh}^{2}\E\LL_{\dimh}(\upsilonvs_{\dimh})\in\R^{\dimtotal\times\dimtotal}\), \(\VFc^{2}=\Cov[\nabla_{\dimp+\dimh}\LL_{\dimh}(\upsilonvs_{\dimh})]\in\R^{\dimtotal\times\dimtotal}\) and \(\upsilonvd=\upsilonvs_{\dimh}\in\R^{\dimtotal}\). Choose \(\rups(\xx)>0\) such that \(\P(\tilde{\upsilonv}_{\dimh},\tilde{\upsilonv}_{\thetavs_{\dimh},\dimh}\in\Ups_{0,\dimh}(\rups(\xx)))\ge 1-\ex^{-\xx}\). 
Then it holds for any \(\dimh\ge\dimh_0\) with probability greater \(1-6\ex^{-\xxn}\)
\begin{EQA}
\bigl\| 
        \DPr_{\dimh} \bigl( \tilde{\thetav}_{\dimh} - \thetavs \bigr) 
        - \xivr_{\dimh}(\upsilonvs_{\dimh}) 
    \bigr\|
    &\le& 
    \,\Excgr(\rups,\xx)+ \alpha(\dimh),
\end{EQA}    
where
 \begin{EQA}   
    \xivr_{\dimh}(\upsilonvs_{\dimh})&\eqdef& \DPr_{\dimh}^{-1}(\nabla_{\thetav}- \A_{\dimh} \HH_{\dimh}^{-1}\nabla_{\etav})\LL(\upsilonvs_{\dimh}).
\end{EQA}
\end{corollary}

For the bias in the Wilks result a bit more work is needed. We can show the following:

\begin{theorem}
\label{theo: wilks bias}
Assume the same as in Theorem~\ref{cor: semi sieve bias}. Pick a radius \(\rups>0\) such that
\begin{EQA}[c]
\P\left(\tilde{\upsilonv}_{\dimh},\tilde{\upsilonv}_{\thetavs_{\dimh},\dimh}\in{\Ups}_{0,\dimh}(\rups)\right)>1-\ex^{-\xx},
\end{EQA}
and set
\begin{EQA}
	\deltar_{\dimp}^{*}
	&\eqdef&
	\left(2+\frac{1+\corrDF}{1-\corrDF}\right)
	\Excgr(\rups,\xx)+\frac{1+\corrDF}{1-\corrDF} \alpha(\dimh).
\end{EQA}
Then we get with probability greater \(1-6\ex^{-\xx}\)
\begin{EQA}
	&& \nquad
	\bigl| 
		2 \Lr(\tilde{\thetav}_{\dimh},\thetavs) - \| \xivr_{\dimh}(\upsilonvs_{\dimh}) \|^{2}
	\bigr|
	\\
	&\le& 
	2\left\{ \rr_{\dimp} + \Excgr\bigl(\rups + \deltar_{\dimp}^{*},\xx\bigr)\right\} 
	\Excgr\bigl(\rups+\deltar_{\dimp}^{*},\xx\bigr)
	\\
	&&
	+ \, \left\{ \deltar_{\dimp}^{*}+\Excgr\bigl(\rups+\deltar_{\dimp}^{*},\xx\bigr)\right\}^{2}
	+ 2 \alpha^{2}(\dimh) + 2\alpha(\dimh)\zz(\xx,\BBr).
\end{EQA}
\end{theorem}
\begin{remark}
Remember that
\begin{EQA}[c]
	\Lr(\thetav)\eqdef \max_{\etav\in\R^{\dimh}}\LL_{\dimh}(\thetav,\etav),
\end{EQA}
where it is important to note that the maximization is restricted to the finite dimensional space \(\R^{\dimh}\).
\end{remark}

\begin{remark}
With Lemma \ref{lem: large dev refinement} we can replace in the above result \(\rups>0\) by \(\rupf\le \rups\) from \eqref{eq: def or refined radius}.
\end{remark}

Now we want to show how this approach allows to prove efficiency of the sieve profile MLE \(\tilde\thetav_{\dimh}\in\R^{\dimp}\). 
From this point we focus on the correctly specified i.i.d. model in which \( \nsize \) denotes the 
sample size and \( \lkh(\thetav,\etav) \) is the log-likelihood for a single observation.
It is known from the convolution theorem (\cite{van1996weak}, Theorem 3.11.2 p. 414, setting \(\kappa(\P_{\upsilonv})=\thetav\)) that the asymptotically optimal variance for regular estimators is given by the inverse of the partial information matrix
\begin{EQA}
    \iif_{\thetav,\etav}
    &=&\left(\Proj\Cov \bigl\{ 
        \nabla \lkh(\upsilonvs)   \bigr\}^{-1}\Proj^{\T}\right)^{-1},
\end{EQA} 
where \( \Proj \) is the orthogonal projection onto the \(\thetav\)-components, and \( \Proj^{\T} \) its adjoint operator. In the case of correct specification we have that \(\DPr^{2}=n \iif_{\thetav,\etav}\).

As the efficient covariance is derived for the score evaluated in the true full target \(\upsilonvs\in l^{2}\) we need further assumptions on the bias:

\begin{description}
\item[\(\bb{(bias')}\)] As \(\dimh\to \infty\) with \(\|\cdot\|\) denoting the spectral norm
\begin{EQA}
\| I-\DPrp(\upsilonvs)^{-1}\DPr(\upsilonvs)^{2}\DPrp(\upsilonvs)^{-1}\|&=& o(1),\\
\| I-\DPrp(\upsilonvs_{\dimh})^{-1}\DPrp(\upsilonvs)^{2}\DPrp(\upsilonvs_{\dimh})^{-1}\|&=& o(1).
 \end{EQA}
\end{description}

\begin{remark}
This paper focuses on the result \ref{theo: main theo finite dim} and thus we do not elaborate on approximation theory. But \cite{AAbias2014} presents conditions on the structure of \(\DF: l^{2}\to l^{2}\) and on the sequence \(\etavs\in l^{2}\) that yield \(\bb{(bias')}\) in an adequate way. 
\end{remark}

Further we need convergence of the covariance of the weighted score.  For this define \(\VPr^{2}_{\dimh,\DF}(\upsilonvs_{\dimh})\eqdef \Cov\left(\score_{\thetav} \lkh_{1}(\upsilonvs_{\dimh})- \A_{\dimh}\HH_{\dimh}^{-2}\score_{\etav} \lkh_{1}(\upsilonvs_{\dimh})\right)\).

\begin{description}
\item[\(\bb{(bias'')}\)] As \(\dimh\to \infty\) with \(\|\cdot\|\) denoting the spectral norm
\begin{EQA}[c]
 \| \DPr_{\dimh}^{-1}\VPr^{2}_{\dimh,\DF}(\upsilonvs_{\dimh})\DPr_{\dimh}^{-1} - I_{\dimp}\|\to 0.
 \end{EQA}
\end{description}

\begin{remark}
This is a condition on how the covariance operator of \(\score_{\dimp+\dimh} \LL(\upsilonv)\in \R^{\dimp+\dimh}\) is affected when it is evaluated in \(\upsilonvs_{\dimh}\in \R^{\dimp+\dimh}\) instead of \(\upsilonvs\in l^{2}\). In the single index example we get \(\bb{(bias'')}\) due to the smoothness of the score.
\end{remark}

 Theorems~\ref{cor: semi sieve bias} and \ref{theo: wilks bias} allow to derive the following corollaries which yield the asymptotic efficiency of \( \tilde{\thetav}_{\dimh} \) and the classical Wilks phenomenon.

\begin{corollary}
\label{cor: efficiency in iid case}
Assume that we have iid observations from \( \P = \P_{\thetavs,\fvs} \) and that the conditions of Theorem~\ref{cor: semi sieve bias} are satisfied. Further let \(\bb{(bias')}\) and \(\bb{(bias'')}\) be satisfied. Assume that there exist sequences \((\dimh_\nsize)\subset\N\) and \(\xxn\to\infty\) with
\begin{EQA}[c]
\label{eq: conditions on approx error in semiefficient}
	\Excgr\bigl(\rups,\xxn\bigr)+ \alpha(\dimh)\to 0,
\end{EQA}
as \( \nsize \to \infty \), where \(\rups(\xxn)\) is chosen such that \(\P(\tilde{\upsilonv}_n,\tilde{\upsilonv}_{\thetavs_{\dimh},\dimh}\in\Ups_{0,\dimh}(\rups(\xxn)))\ge 1-\ex^{-\xxn}\). Then as \( \nsize \to \infty \)
\begin{EQA}[ccl]
    (\nsize \iif_{\thetav,\fv})^{1/2} \bigl( \tilde{\thetav}_{\dimh} - \thetavs \bigr)
    - \xivr
    & \toP &
    0,
    \\
    (\nsize \iif_{\thetav,\fv})^{1/2} \bigl( \tilde{\thetav}_{\dimh} - \thetavs \bigr)
    & \tow &
    \ND(0,\Id_{\dimp}).
\end{EQA}
\end{corollary}

\begin{corollary}
\label{cor: wilks in iid case}
Assume that we have iid observations from \( \P = \P_{\thetavs,\fvs} \) and that the conditions of Theorem~\ref{cor: semi sieve bias} are satisfied.  Further let \(\bb{(bias')}\) be satisfied. Assume there exist a sequences \((\dimh_\nsize)\subset\N\) and \(\xxn\to\infty\) with 
\begin{EQA}[c]
\label{eq: conditions on approx error in semiefficientiid}
	\Excgr(\rups,\xx) \sqrt{\xxn+\dimp}
	+ \alpha^{2}(\dimh) 
	+ \alpha(\dimh)\sqrt{\xxn+\dimp}
	\to 0,
\end{EQA}
as \( \nsize \to \infty \), where \(\rups(\xxn)\) is chosen such that \(\P(\tilde{\upsilonv}_n,\tilde{\upsilonv}_{\thetavs_{\dimh},\dimh}\in\Ups_{0,\dimh}(\rups(\xxn)))\ge 1-\ex^{-\xxn}\). Then as \( \nsize \to \infty \)
\begin{EQA}[ccl]
   2 \Lr(\tilde{\thetav}_{\dimh},\thetavs)
    & \tow &
    \chi^{2}_{\dimp}.
\end{EQA}
\end{corollary}

\subsection{A single-index model}
\label{sec: application to single index model}
We illustrate how the results from Section~\ref{sec: main 
results} and the last statement can be derived for single-index modeling. 
Consider the following model
\begin{EQA}[c]
\label{eq: single index model introduced}
    \Yv_{i}
    =
    \fs(\Xv_{i}^{\T} \thetavs) + \varepsilon_{i}, 
    \qquad 
    i=1,...,\nsize,
\end{EQA}
or some \(\fs:\R\to \R\) and \(\thetavs\in S_{1}^{\dimp,+}\subset\R^{\dimp}\), i.i.d errors \(\varepsilon_i\in\R\) and \(\Var(\varepsilon_i)=\sigma^{2}\) and i.i.d random variables \(\Xv_{i}\in \R^{\dimp}\) with distribution denoted by \(\P^{\Xv}\). The single-index model is widely applied in statistics. 
For example in econometric studies it serves as a compromise between too restrictive 
parametric models and flexible but hardly estimable purely nonparametric models. 
Usually the statistical inference focuses on estimating the 
index vector \( \thetavs \). 
A lot of research has already been done in this field. 
For instance, \cite{DelecroixHaerdleHristache} show the asymptotic efficiency of the 
general semiparametric maximum-likelihood estimator for particular examples and in 
\cite{HaerdleHallIchimura} the right choice of bandwidth for the nonparametric 
estimation of the link function is analyzed.

To ensure identifiability of \(\thetavs\in\R^{\dimp}\) we assume that it lies in the half sphere \(S_{1}^{\dimp,+}\eqdef\{\thetav\in\R^{\dimp}:\, \|\thetav\|=1,\, \theta_1> 0\}\subset\R^{\dimp}\). For simplicity we assume that the support of the \(\Xv_{i}\in \R^{\dimp}\) is contained in the ball of radius \(s_{\Xv}>0\). 
Further we assume that \(\fs\in \{f:[-s_{\Xv},s_{\Xv}]\mapsto \R\} \) can be well approximated by an orthonormal \(C^{2}\)-Daubechies-wavelet basis, i.e. for a suitable function \(\basX_{0}\eqdef\psi:[-s_{\Xv},s_{\Xv}]\mapsto \R\) we set for \(k=r_k 12+2^{r_k}+j_k\) with \(r_k\in\N_0\) and \(j_k\in\{0,\ldots,11+2^{r_k}\}\)
\begin{EQA}[c]
\basX_{k}(t)=2^{r_k/2}\psi\left(2^{r_k}(t-2j_ks_{\Xv})\right),\, k\in\N.
\end{EQA}

Our aim is to analyze the properties of the profile MLE 
\begin{EQA}[c]
    \tilde{\thetav}_{\dimh}
    \eqdef
    \argmax_{\thetav} \max_{\etav \in \R^{\dimh}} \LL_{\dimh}(\thetav,\etav) ,
\label{ttSI}
\end{EQA}    
where
\begin{EQA}[c]
    \LL_{\dimh}(\thetav,\etav)
    =
    - \frac{1}{2} 
    \sum_{i=1}^{\nsize} \Bigl| \Yv_{i} - \sum_{k=0}^{m} \etav_{k} \basX_{k}(\Xv_{i}^{\T} \thetav) \Bigr|^{2}.
\end{EQA}
\cite{Ichimura} analyzed a very similar estimator in a more general setting based on 
a kernel estimation of \( \E \bigl[ \Yv \cond \fs(\thetav^{\T} \Xv) \bigr] \) 
instead of using a parametric sieve approximation \( \sum_{k=0}^{m} \etav_{k} \basX_{k} \). 
He showed \( \sqrt{\nsize} \)-consistency and asymptotic normality of the proposed 
estimator.

To apply the technique presented in the previous section we need a list of assumptions. We denote this list of conditions by \((\mathbf{\mathcal A})\). We start with conditions on the regressors \(\Xv\in\R^{\dimp}\):

\begin{description} 
 \item[\((\mathbf{Cond}_{\Xv})\)]
  The measure \(\P^{\Xv}\) is absolutely continuous with respect to the Lebesgue measure. The Lebesgue density \(d_{\Xv}: \R^{\dimp}\to \R\) of \(\P^{\Xv}\) is only positive on the ball \(B_{s_{\xv}}(0)\subset \R^{\dimp}\) and Lipschitz continuous with Lipschitz constant \(L_{d_{\Xv}}>0\). Further we assume that for any \(\thetav\perp\thetavs\) with \(\|\thetav\|=1\) we have \(\Var\Big(\Xv^\T\thetav\Big|\Xv^\T\thetavs\Big)>\sigma^2_{\Xv|\thetavs}\) for some constant \(\sigma^2_{\Xv|\thetavs}>0\) that does not depend on \(\Xv^\T\thetavs\in\R\). Also assume that that the density \(d_{\Xv}: \R^{\dimp}\to \R\) of the regressors satisfies \(d_{\Xv}>c_{d_{\Xv}}>0\) on  \(B_{s_{\xv}}(0)\subset \R^{\dimp}\)
\end{description}

\begin{remark}
We only assume bounded support of the regressors \((\Xv_i)\subset \R^{\dimp}\) for simplicity. This condition could be relaxed to a qualified probabilistic deviation bound of the kind \(\P(\|\Xv\|\ge s_{\Xv}+\xx)\le \ex^{-\xx}\). Further \(\Var\Big(\Xv^\T\thetavd\Big|\Xv^\T\thetavs\Big)=0\) would mean that \(\Xv^\T\thetavd=h(\Xv^\T\thetavs)\) for some function \(h:\R\to\R\). But then we would have for any \((\alpha,\beta)\in \R^2\) with \(\alpha^2+\beta^2=1\) that
\begin{EQA}[c]
f(\Xv^\T(\alpha\thetavs+\beta\thetavd))=f(\alpha\Xv^\T\thetavs+\beta h(\Xv^\T\thetavs))\eqdef g_{\alpha,\beta}(\Xv^\T\thetavs),
\end{EQA}
such that the problem would no longer be identify-able. We bound \(d_{\Xv}>c_{d_{\Xv}}>0\) on  \(B_{s_{\xv}}(0)\subset \R^{\dimp}\) as a simple way to ensure identifyability. 
\end{remark}

Of course we need some regularity of the link function \(\fs\in \{f:[-s_{\Xv},s_{\Xv}]\mapsto \R\} \):
\begin{description} 
  
  \item[\( (\mathbf{Cond}_{\fs}) \)]
   For some \(\etavs\in l^2\)
   \begin{EQA}[c]
   \label{eq: expansion of indexfunciton}
    \fs=\fs_{\etavs}=\sum_{k=1}^\infty \eta^*_k\basX_{k},
   \end{EQA}
   where with some \(\alpha>2\) and a constant \(C_{\|\etavs\|}>0\)
   \begin{EQA}[c]\label{eq: smoothness of fs}
    \sum_{l=0}^\infty l^{2\alpha}{\eta^*_{l}}^2 \le C_{\|\etavs\|}^2< \infty.
   \end{EQA} 
\end{description}

For the large deviations of the MLE we need the following condition:
\begin{description} 
 \item[\((\mathbf{Cond}_{\Xv\thetavs})\)]
  It holds true that \(\P(|\fs_{\etavs}'(\Xv^\T\thetavs)|> c_{\fs_{\etavs}'})>c_{\P\fs'}\) for some \(c_{\fs_{\etavs}'},c_{\P\fs'}>0\). 
\end{description}

\begin{remark}
 Note that a condition of this kind is necessary to ensure identify ability. Otherwise the function \(\fs\) would be \(\P^{\Xv}\)-almost surely constant. But for a constant function \(C(x)\equiv c\) any \(\thetav\in\R^{\dimp}\) solves \(C(\Xv^\T \thetav)= c\).
\end{remark}

To be able to apply the finite sample device we need constraints on the moments of the additive noise:

\begin{description} 
 \item[\((\mathbf{Cond}_{\varepsilon})\)] 
  The errors \((\varepsilon_{i})\in \R\) are i.i.d. with \(\E[\varepsilon_{i}]=0\), \(\Cov(\varepsilon_{i})=\sigma^2\) and satisfy for all \(|\mu|\le\tilde g\) for some \(\tilde g>0\) and some \(\tilde \nu_{\rr}>0\)
  \begin{EQA}[c]
    \log\E[\exp\left\{ \mu \varepsilon_{1} \right\}]\le \tilde \nu_{\rr}^{2}\mu^{2}/2.
  \end{EQA}
\end{description}

Finally to be able to control the large deviations of the MLE we impose
 \begin{description}  
    \item[\((\mathbf{Cond}_{\Ups})\)]  \(\Ups\subseteq\Upss(\sqrt n\rr^{\circ})\subset \R^{\dimp+\dimh}\) with \(\rr^{\circ}\in\R\). Such that \( d_{\Ups}\eqdef diam(\Ups)<\infty\).
\end{description}

If these conditions denoted by \((\mathbf{\mathcal A})\) are met we can proof the following result:
\begin{proposition}
\label{prop: semi sieve bias}
Assume \((\mathbf{\mathcal A})\). If \(\alpha= 2+\eps\), if \({\dimtotal}^{5}/ n\to 0\) but \({\dimtotal}^{5+2\eps}/ n\to \infty\) for some \(\eps>0\)  and if \(n\in\N\) is large enough
we get 
\begin{EQA}
	\Excgr(\rups,\xx)
	\lesssim 
	\frac{(\dimtotal + \xx)^{5/2}}{\sqrt{\nsize}},\quad
	&
	\alpha(\dimh) 
	\lesssim 
	\sqrt{\nsize}\dimh^{-\alpha-1/2},\quad
	& \beta(\dimh)\lesssim \dimh^{-1}.
\end{EQA}
Further as \( \nsize \to \infty \)
\begin{EQA}[rcl]
    \bigl\| 
        \DPr \bigl( \tilde{\thetav}_{\dimh} - \thetavs \bigr) 
        - \xivr_{\dimh}(\upsilonvs_{\dimh}) 
    \bigr\|
    & \toP &
    0,
    \\
    \DPr \bigl( \tilde{\thetav}_{\dimh} - \thetavs \bigr)
    & \tow &
    \ND(0,\sigma^{2}\Id_{\dimp}),\\
      2 \Lr(\tilde{\thetav}_{\dimh},\thetavs)
  & \tow &
  \chi^{2}_{\dimp}.
\end{EQA}
\end{proposition}
For details see \cite{Andresensingleindex}.

\appendix
\section{Appendix}

\section{Deviation bounds for quadratic forms}
\label{ap: deviation of quadratic forms}
The following general result from \cite{SP2011} helps
to control the deviation for quadratic forms of type \( \| \BB \xiv \|^{2} \)
for a given positive matrix \( \BB \) and a random vector \( \xiv \). 
It will be used several times in our proofs.
Suppose that 
\begin{EQA}[c]
    \log \E \exp\bigl( \gammav^{\T} \xiv \bigr) 
    \le 
    \| \gammav \|^{2}/2,
    \qquad 
    \gammav \in \R^{\dimp}, \, \| \gammav \| \le \gm .
\label{expgamgm} 
\end{EQA}
For a symmetric matrix \( \BB \), define 
\begin{EQA}[c]
    \dimA = \tr (\BB^{2}) , 
    \qquad 
    \vA^{2} = 2 \tr(\BB^{4}),
    \qquad 
    \lambdaB \eqdef \| \BB^{2} \|_{\infty} \eqdef \lambda_{\max}(\BB^{2}) .
\label{dimAvAlb}
\end{EQA}
For ease of presentation, suppose that \( \gm^{2} \ge 2 \dimB \).
The other case only changes the constants in the inequalities. 
Note that \( \| \xiv \|^{2} = \etav^{\T} \BB \, \etav \).
Define \( \muc = 2/3 \) and
\begin{EQA}
    \gmc
    & \eqdef &
    \sqrt{\gm^{2} - \muc \dimB} ,
    \\
    2(\xxc+2)
    & \eqdef &
    (\gm^{2}/\muc - \dimB)/\lambdaB + \log \det \bigl( \Id_{\dimp} - \muc \BB/\lambdaB \bigr) .
\label{yycgmcxxcAd}
\end{EQA}

\begin{proposition}
\label{LLbrevelocro}   
\label{theo: dev bounds quad forms}
Let \( (E\!D_{0}) \) hold with \( \nunu = 1 \) and 
\( \gmb^{2} \ge 2 \dimB \).
Then
for each \( \xx > 0 \)
\begin{EQA}
    \P\bigl( \| \xiv \| \ge \zz(\xx,\BB) \bigr)
    & \le &
    2 \ex^{-\xx} ,
\label{PxivbzzBBro}
\end{EQA}    
where \( \zz(\xx,\BB) \) is defined by
\begin{EQA}
\label{PzzxxpBro}
    \zz^{2}(\BB,\xx)
    & \eqdef &
    \begin{cases}
        \dimB + 2 \vA_{\BB} (\xx+1)^{1/2}, &  \xx+1 \le \vA_{\BB}/(18 \lambdaB) , \\
        \dimB + 6 \lambdaB (\xx+1), & \vA_{\BB}/(18 \lambdaB) < \xx+1 \le \xxc+2 , \\
        \bigl| \yyc + 2 \lambdaB (\xx - \xxc + 1)/\gmc \bigr|^{2}, & \xx > \xxc+1 ,
    \end{cases}
\label{zzxxppdBlro}
\end{EQA}    
with \( \yyc^{2} \leq \dimB + 6 \lambdaB (\xxc+2) \).
\end{proposition}

Depending on the value \( \xx \), we observe three types of tail behavior of the 
quadratic form \( \| \xiv \|^{2} \). 
The sub-Gaussian regime for 
\( \xx+1 \le \vA_{\BB}/(18 \lambdaB) \) 
and the Poissonian regime for \( \xx \le \xxc + 1 \) 
are similar to the case of a Gaussian quadratic form.
The value \( \xxc \) from \eqref{yycgmcxxcAd} is of order \( \gm^{2} \).
In all our results we suppose that \( \gm^{2} \) and hence, \( \xxc \) 
is sufficiently large and 
the quadratic form \( \| \xiv \|^{2} \) can be bounded with a dominating 
probability by \( \dimB + 6 \lambdaB (\xx+1) \) for a proper \( \xx \).
%
We refer to \cite{SP2011} for the proof of this and related results, further 
discussion and references.


\section{Proofs}
\label{Sproofsa}
This section collects the proofs of the results in chronological order.

\subsection{Proof of Theorem~\ref{theo: main theo finite dim}}
We proof the theorem in three steps. First we present a local approximation result for the full parameter \(\upsilonv \in \Ups\), then we proof \eqref{eq: Fisher in main theo} and \eqref{eq: Wilks in main theo} on a suitable subset \(\LCS(\rr)\subset \Omega\) of the underlying probability space and finally we show that \(\LCS(\rr)\subset \Omega\) is of dominating probability.

\subsubsection{Local approximations in the full space}
\label{sec: local approx full space}
We start the proof with an auxiliary theorems that is of interest on its own. Define the approximation error
\begin{EQA}[c]
    \Excgr(\rr,\xx)
    \eqdef
    \left(\rddelta(\rr) + \, 6\nuno \rhor \zzQ(\xx,\dimtotal) \right)\rr,
\label{Exceqrrrhdef}
\end{EQA}
and the normalized stochastic gradient gap
\begin{EQA}[c]
\UU(\upsilonv)=\DF^{-1}\Big(\nabla\zetav(\upsilonv) - \nabla\zetav(\upsilonvs)\Big).
\end{EQA}

\begin{theorem}
 \label{theo: local expansion full parameter}
Assume that the condition \( \bb{(\LL_{0})} \) is fulfilled (where \(\dimtotal=\dimh+\dimp\)). Then on the set 
\begin{EQA}
\label{eq: assumpion on bound for norm of sup of stoch grad}
	\mathcal M(\xx)
	&\eqdef& 
	\{\tilde{\upsilonv},\tilde{\upsilonv}_{\thetavs},\DF^{-2}\nabla \LL(\upsilonvs)+\upsilonvs
		\in \Upss(\rr)\}
	\cap\{\|\DF^{-1}\nabla \LL(\upsilonvs)\|\le \zz(\xx,\BB)\}
	\\
	&&
	\cap\left\{\sup_{\upsilonv\in\Upss(\rr)}\|\UU(\upsilonv)\|
	\le 
	6\nuno \rhor \zzQ(\xx,\dimtotal) \right\},
\end{EQA}
we have for \(\rr \ge \zz(\xx,\BB)\)
\begin{EQA}[c]
\|\DF(\tilde{\upsilonv}-\upsilonvs)-\xiv\|\le \Excgr(\rr,\xx).
\end{EQA}
\end{theorem}

\begin{proof}
Define
\begin{EQA}[c]
	\alp(\upsilonv,\upsilonvs)
	\eqdef
	\LL(\upsilonv) - \LL(\upsilonvs) 
    - (\upsilonv - \upsilonvs)^{\T} \nabla \LL(\upsilonvs) 
    + \frac{1}{2} \| \DF (\upsilonv - \upsilonvs) \|^{2},
\label{eq: error in taylor approx}
\end{EQA}
and note that by definition
\begin{EQA}
\LL(\upsilonv,\upsilonvs)&=&\nabla \LL(\upsilonvs)(\upsilonv-\upsilonvs)-\|\DF(\upsilonv-\upsilonvs)\|^{2}/2+\alp(\upsilonv,\upsilonvs).
\end{EQA}
Setting \(\nabla\LL(\tilde{\upsilonv})=0\) we find that \(\tilde{\upsilonv} \) satisfies
\begin{EQA}[c]
\DF(\tilde{\upsilonv}-\upsilonvs)-\DF^{-1}\nabla \LL(\upsilonvs)=\DF^{-1}\nabla \alp(\tilde{\upsilonv},\upsilonvs).
\end{EQA}
This gives
\begin{EQA}[c]
\|\DF(\tilde{\upsilonv}-\upsilonvs)-\DF^{-1}\nabla \LL(\upsilonvs)\|^{2}= (\DF(\tilde{\upsilonv}-\upsilonvs)-\DF^{-1}\nabla \LL(\upsilonvs))^{\T}\DF^{-1}\nabla \alp(\tilde{\upsilonv},\upsilonvs).
\end{EQA}
As we have on \(\mathcal M(\xx)\subset \Omega\) that \(\tilde{\upsilonv},\DF^{-2}\nabla \LL(\upsilonvs)+\upsilonvs \in\Upss(\rr)\) 
\begin{EQA}[c]
	\|\DF(\tilde{\upsilonv}-\upsilonvs)-\DF^{-1}\nabla \LL(\upsilonvs)\|
	\le 
	\|\DF^{-1}\nabla \alp(\tilde{\upsilonv},\upsilonvs)\|,
\end{EQA}
So it suffices to show that 
\begin{EQA}[c]
	\sup_{\upsilonv\in\Upss(\rr)} \|\DF^{-1} \nabla \alp( \upsilonv,\upsilonvs)\|
	\le \Excgr(\rr,\xx),
\end{EQA}
where
\begin{EQA}
\label{eq: def of U process}
	\DF^{-1}\nabla \alp( \upsilonv,\upsilonvs)
	& =&
	\DF^{-1} \bigl\{ 
		\nabla \LL(\upsilonv) - \nabla \LL(\upsilonvs) - \DF^{2} \, (\upsilonv - \upsilonvs) 
	\bigr\}.
\label{UPupsnm}
\end{EQA}
Note that by condition \( \bb{(\LL_{0})} \) and Taylor expansion
\begin{EQA}
\label{eq: bound for expected value of approx error}
	&&\nquad
	\sup_{\upsilonv\in\Upss(\rr)} \|\E\DF^{-1}\nabla \alp( \upsilonv,\upsilonvs)\|
	\\
  	&=& 
	\sup_{\upsilonv\in\Upss(\rr)}\|\DF^{-1}\nabla \E\LL(\upsilonv) 
	- \DF^{-1}\nabla \E\LL(\upsilonvs)-\DF \, (\upsilonv - \upsilonvs)\|
	\\
  	&\le& 
	\sup_{\upsilonv\in\Upss(\rr)}\|\DF^{-1} 
		\nabla^{2} \E\LL(\upsilonv)^{2}\DF^{-1} - \Id_{\dimtotal}\| \rr
	\\
  	&\le& 
	\delta(\rr)\rr.
\end{EQA}
For the remainder we use assumption \eqref{eq: assumpion on bound for norm of sup of stoch grad}. This gives the claim.
\end{proof}


\subsubsection{Proof of Theorem~\ref{theo: main theo finite dim} on the set \(\LCS(\rups,\xx)\)}
We now proof Theorem~\ref{theo: main theo finite dim} on a suitable subset \(\LCS(\rups,\xx)\subset \Omega\). For this purpose fix some radius \(\rups(\xx)>0\) that ensures dominating probability for the event \(\{\tilde{\upsilonv},\tilde{\upsilonv}_{\thetavs}\in\Upss(\rups)\}\). Define \( \LCS(\rups,\xx) \subseteq \Omega\) as 
\begin{EQA}
\label{eq: def of local set}
    \LCS(\rups,\xx)
    \eqdef& \bigl\{ &
    \tilde{\upsilonv},\tilde{\upsilonv}_{\thetavs}\in \Upss(\rups),
    \\
    &&    \| \DF^{-1} \score \| \le \zz(\xx,\BB),\, 
        \| \HH^{-1} \score_{\etav} \| \le \zz(\xx,\BB)/(1-\corrDF),\, \|\DPr^{-1}\scorer\|\le \zz(\xx,\BBr) 
    \bigr\}\\	
    &\cap&\bigg\{\sup_{\upsilonv\in\Upss(\rr)}\|\UU(\upsilonv)\|\le 6\nuno \rhor \zzQ(\xx,\dimtotal) \bigg\},
\end{EQA}
where
\begin{EQA}[c]
\tilde{\upsilonv}_{\thetavs}\eqdef \argmax_{\Proj\upsilonv=\thetavs}\LL(\upsilonv).
\end{EQA}

In the following we will derive statements that hold true on this set \( \LCS(\rups,\xx) \subseteq \Omega\) which is of probability greater \(1-6\ex^{-\xx}\) as we show later. Note that on \( \LCS(\rups,\xx) \subseteq \Omega\) we have \( \hat{\upsilonv}+\upsilonvs\in \Upss(\zz(\xx,\BB))\) and \((\thetavs,\hat{\etav}+\etavs)\in \Upss(\zz(\xx,\BB)/(1-\corrDF))\), where
\begin{EQA}
    \hat{\upsilonv} 
    & \eqdef &
    \DF^{-2} \score
    =
    \argmin_{\upsilonv} \La(\upsilonv,\upsilonvs)-\upsilonvs,
    \\
    \hat{\etav} 
    & \eqdef &
    \HH^{-2} \score_{\etav}
    =
    \argmin_{\etav} \La((\thetavs,\etav),\upsilonvs)-\etavs .
\label{tupsrdm}
\end{EQA}
Remember the notation \( \score \eqdef \nabla \LL(\upsilonvs) \) and the definition of \( \scorer_{\thetav} \) and \( \DPr^{2} \):
\begin{EQA}
     \scorer_{\thetav} 
     & \eqdef &
     \score_{\thetav} - \A \HH^{-2} \score_{\etav},
     \\
     \DPr^{2}
     & \eqdef &
     \DP^{2} - \A \HH^{-2} \A^{\T} .
\label{DPrbAbHHbappp}
\end{EQA}

To prove the finite sample Fisher expansion \eqref{eq: Fisher in main theo} we take the expansion for the whole parameter vector \( \tilde{\upsilonv} \in \Ups \) from Theorem~\ref{theo: local expansion full parameter}, i.e. that on the set \(\LCS(\rups,\xx)\subset \Omega\) we have 
\begin{EQA}[c]
\label{eq: full expansion}
\|\DF(\tilde{\upsilonv}-\upsilonvs)-\DF^{-1}\score\|\le \Excgr(\rups,\xx).
\end{EQA}
It remains to note that for any 
\( \uv \in \R^{\dimp} \), \( \etav \in \R^{\dimh} \), and 
\( \wv = (\uv,\etav) \in \R^{\dimtotal} \), it holds
with \( \gammav \eqdef \etav + \HH^{-2} \A^{\T} \uv \in \R^{\dimh} \)
\begin{EQA}[c]
    \bigl\| \DF \wv \|^{2}
    = 
    \bigl\| \DPr \uv \bigr\|^{2} + \bigl\| \HH \gammav \bigr\|^{2}
    \ge 
    \bigl\| \DPr \uv \bigr\|^{2} .
\label{DFcDPcuvse}
\end{EQA}    
Further \( \Proj \DF^{-2} \score = \DPr^{-2} \scorer \). 
This implies for \( \wv = \tilde{\upsilonv} - \upsilonvs \) by \eqref{DFcDPcuvse}
\begin{EQA}[c]
    \bigl\| \DPr (\tilde{\thetav} - \thetavs) - \DPr^{-1} \scorer \bigr\| 
    =
    \bigl\| \DPr (\tilde{\thetav} - \thetavs - \DPr^{-2} \scorer) \bigr\| 
    \le 
    \bigl\| \DF (\wv - \DF^{-2} \score) \bigr\|,
\label{DPbtcse}
\end{EQA}    
and the assertion follows on the set \( \LCS(\rups,\xx) \).

Before we show that the bound \eqref{eq: Wilks in main theo} is fulfilled 
on the set \( \LCS(\rups,\xx) \) from \eqref{eq: def of local set} we present a list of auxiliary Lemmas.
\begin{lemma}
\label{LLabmProjse}
It holds on the set 
\( \bigl\{ \| \DF^{-1} \score \| \le \rupf, \, 
\| \HH^{-1} \score_{\etav} \| \le  \rupf \bigr\} \)
\begin{EQA}
    \sup_{\upsilonv} \La(\upsilonv,\upsilonvs)
    &=&
    \sup_{\upsilonv \in \Upss(\rupf)} \La(\upsilonv,\upsilonvs)
    =
    \frac{1}{2} \bigl\| \DF^{-1} \score \bigr\|^{2} ,\label{supLabapp}
    \\
    \sup_{\upsilonv \in \Upss(\rupf): \, \Proj\upsilonv=\thetavs} \La(\upsilonv,\upsilonvs)
    &=&
    \sup_{\upsilonv: \, \Proj\upsilonv=\thetavs} \La(\upsilonv,\upsilonvs) 
    =
    \frac{1}{2} \bigl\| \HH^{-1} \score_{\etav} \bigr\|^{2},\label{supprojLabapp0}
    \\
    \sup_{\upsilonv} \La(\upsilonv,\upsilonvs)
    - \sup_{\upsilonv: \, \Proj\upsilonv=\thetavs} \La(\upsilonv,\upsilonvs)
    &=&
    \frac{1}{2} \bigl\| \DPr^{-1} \scorer_{\thetav} \bigr\|^{2}.
\label{supLase}
\end{EQA}    
\end{lemma}

\begin{proof}
First consider the adaptive cases with \( \A = 0 \) yielding \( \DPr^{2} = \DP^{2} \)
and \( \scorer_{\thetav} = \score_{\thetav} \).
Then the process \( \La(\upsilonv,\upsilonvs) \) can be decomposed as
\begin{EQA}
    \La(\upsilonv,\upsilonvs)
    &=&
    (\thetav - \thetavs)^{\T} \score_{\thetav} 
    - \frac{1}{2} \| \DP (\thetav - \thetavs) \|^{2}
    \\
    &+&
    (\etav - \etavs)^{\T} \score_{\etav} 
    - \frac{1}{2} \| \HH (\etav - \etavs) \|^{2},
\label{Labadaptapp}
\end{EQA}    
and the partial optimization subject to \( \thetav = \thetavs \) yields the result
\eqref{supprojLabapp0}.
Note that the constrained maximum is attained at 
\( \etav = \etavs + \HH^{-2} \score_{\etav} \).

The general case can be reduced to the adaptive one by the change of variable.
With \( \gammav \eqdef \etav - \etavs + \HH^{-2} \A^{\T} (\thetav - \thetavs) \),
one can represent \( \La(\upsilonv,\upsilonvs) \) in the form
\begin{EQA}
    \La(\upsilonv,\upsilonvs)
    &=&
    (\thetav-\thetavs)^{\T} \scorer_{\thetav} 
    - \| \DPr(\thetav - \thetavs) \|^{2}/2 
    + \gammav^{\T} \nabla_{\etav} 
    - \| \HH \gammav \|^{2}/2,
\label{decomposed expression for linear problem}
\end{EQA}
which corresponds to the decomposition in the adaptive case.
\end{proof}

\begin{lemma}
\label{lem: distance mle and cond mle}
We have on the set \( \LCS(\rups,\xx) \) that
\begin{EQA}[c]
\|\DF(\tilde{\upsilonv}-\tilde{\upsilonv}_{\thetavs})\|\le \left(2+\frac{1+\corrDF}{1-\corrDF}\right)\Excgr(\rupf,\xx)+\frac{1+\corrDF}{1-\corrDF}\zz(\xx,\BBr)\eqdef \rr_{\dimp}.
\end{EQA}
\end{lemma}

\begin{proof}
The first order criteria of maximality are satisfied:
\begin{EQA}[c]
\nabla_{\etav}\LL(\tilde{\upsilonv})=\nabla_{\etav}\LL(\tilde{\upsilonv}_{\thetavs})=0.
\end{EQA}
Remember the definition of \(\alp(\upsilonv,\upsilonvs)\) in \eqref{eq: error in taylor approx}
\begin{EQA}[c]
	\alp(\upsilonv,\upsilonvs)
	\eqdef
	\LL(\upsilonv) - \LL(\upsilonvs) 
    - (\upsilonv - \upsilonvs)^{\T} \nabla \LL(\upsilonvs) 
    + \frac{1}{2} \| \DF (\upsilonv - \upsilonvs) \|^{2} .
\end{EQA}
We find
\begin{EQA}
0&=&\HH^{2}(\tilde\etav-\etavs)-2\A^{\T}(\tilde\thetav-\thetavs)-\nabla_{\etav} \LL(\upsilonvs)+\nabla_{\etav} \alp(\tilde{\upsilonv},\upsilonvs)\\
&=&\HH^{2}(\tilde\etav_{\thetavs}-\etavs)-\nabla_{\etav} \LL(\upsilonvs)+\nabla_{\etav} \alp(\tilde{\upsilonv}_{\thetavs},\upsilonvs),
\end{EQA}
from which we derive
\begin{EQA}[c]
\HH(\tilde\etav_{\thetavs}-\tilde\etav)= \HH^{-1}\nabla_{\etav} \alp(\tilde{\upsilonv},\upsilonvs)-\HH^{-1}\nabla_{\etav} \alp(\tilde{\upsilonv}_{\thetavs},\upsilonvs)+2\HH^{-1}\A^{\T}(\tilde\thetav-\thetavs),
\end{EQA}
where
\begin{EQA}[c]
\tilde \etav_{\thetavs}\eqdef \argmax_{\substack{\etav\in\R^{\dimh}\\ (\thetavs,\etav)\in\Upsilon}}\LL(\thetavs,\etav).
\end{EQA}
With the same arguments as in the proof of \ref{theo: local expansion full parameter} we infer using that \(\tilde{\upsilonv},\tilde{\upsilonv}_{\thetavs}\in\Upss(\rupf)\)
\begin{EQA}[c]
\label{eq: bound for hh diff eta mle cond and non cond}
\|\HH(\tilde\etav_{\thetavs}-\tilde\etav)\|\le 2\Excgr(\rupf,\xx)+\corrDF\|\DP(\tilde\thetav-\thetavs)\|.
\end{EQA}
This gives on \( \LCS(\rups,\xx) \)
\begin{EQA}
\|\DF(\tilde{\upsilonv}-\tilde{\upsilonv}_{\thetavs})\|&\le& \|\HH(\tilde\etav_{\thetavs}-\tilde\etav)\|+\|\DP(\tilde\thetav-\thetavs)\|\\
	&\le&2\Excgr(\rupf,\xx)+(1+\corrDF)\|\DP(\tilde\thetav-\thetavs)\|\\
	&\le&2\Excgr(\rupf,\xx)+(1+\corrDF)\|\DP\DPr^{-1}\|\|\DPr(\tilde\thetav-\thetavs)\|\\
	&\le&2\Excgr(\rupf,\xx)+\frac{1+\corrDF}{1-\corrDF }\left(\|\DPr^{-1}\scorer\|+\Excgr(\rupf,\xx)\right)\\
	&\le& \left(2+\frac{1+\corrDF}{1-\corrDF} \right)\Excgr(\rupf,\xx)+\frac{1+\corrDF}{1-\corrDF} \zz(\xx,\BBr).
\end{EQA}
This gives the claim. 
\end{proof}

\begin{lemma}
\label{lem: bound on alpha distance}
With \(\alp(\upsilonv,\upsilonvs)\) defined as in \eqref{eq: error in taylor approx} we have on the set \(\LCS(\rups,\xx)\)
\begin{EQA}[c]
\sup_{\substack{(\upsilonv,\upsilonvc)\in\Upss(\rupf)^{2}\\ \|\DF(\upsilonv-\upsilonvc)\|\le \rr_{\dimp}}}\left\{\alp(\upsilonv,\upsilonvs)
    - \alp(\upsilonvc,\upsilonvs) \right\}\le \rr_{\dimp}\Excgr(\rupf,\xx).
\end{EQA}
\end{lemma} 

\begin{proof}
Take any \((\upsilonv,\upsilonvc)\in\Upss(\rupf)^{2}\) with \( \|\DF(\upsilonv-\upsilonvc)\|\le \rr_{\dimp}\). With \(\upsilonvd\in\Upss(\rupf)\) in the convex hull of \(\upsilonv,\upsilonvc\in \Upsilon\) and with \(\UP(\upsilonv)\in\R\) from \eqref{eq: def of U process} we find
\begin{EQA}
|\alp(\upsilonv,\upsilonvs)
    - \alp(\upsilonvc,\upsilonvs)| &=& |(\upsilonv-\upsilonvc)^{\T}\nabla\alp(\upsilonvd,\upsilonvs)|\\
    	&\le& \|\DF(\upsilonv-\upsilonvc)\|\sup_{\upsilonv\in\Upss(\rr)}\|\UP(\upsilonv)\|\le  \rr_{\dimp}\Excgr(\rupf,\xx),
\end{EQA}
where we used the arguments of the proof of Theorem~\ref{theo: local expansion full parameter} to bound \(\UP(\upsilonv)\in\R\). This gives the claim.
\end{proof}

Now we show that the bound \eqref{eq: Wilks in main theo} is fulfilled 
on the set \( \LCS(\rups,\xx) \) from \eqref{eq: def of local set}.
First note that the same arguments which lead to equation \eqref{eq: refined bound for mle} work for \(\tilde{\upsilonv}_{\thetavs}\in\R^{\dimtotal}\) (with a modified version of Theorem~\ref{theo: local expansion full parameter}), which yields on \( \LCS(\rups,\xx)\subseteq \Omega\)
\begin{EQA}[c]
\label{eq: refined bound for thetavs mle}
\|\DF(\tilde{\upsilonv}_{\thetavs}-\upsilonvs)\|\le \zz(\xx,\BB)/(1-\corrDF)+\Excgr(\rups,\xx),
\end{EQA}
such that \( \LCS(\rups,\xx)\subseteq\{ \tilde{\upsilonv}, \tilde{\upsilonv}_{\thetavs} 
\in \Upss(\rupf)\}\subset\Omega \) where \(\rupf\eqdef\zz(\xx,\BB)/(1-\corrDF)+\Excgr(\rups,\xx)\).
\begin{remark}
Note that such an adaptation of Theorem~\ref{theo: local expansion full parameter} works out because the set from \eqref{eq: assumpion on bound for norm of sup of stoch grad} is contained in
\begin{EQA}[c]
\{\tilde{\upsilonv}_{\thetavs}\in\Upss(\rr)\}\cap\left\{\sup_{\substack{\etav\in\R^{\dimh}\\ (\thetavs,\etav)\in\Upss(\rr)}}\|\UU(\upsilonv)\|\le 6\nuno \rhor \zzQ(\xx,\dimtotal) \right\},
\end{EQA}
and because \(\bb{(\LL_0)}\) implies the same for the matrix valued function \(\nabla_{\etav}^{2}\E\LL(\cdot)\).
\end{remark}

We write
\begin{EQA}[c]
    \Lr(\tilde{\thetav},\thetavs)
   =
   \LL(\tilde{\upsilonv},\upsilonvs)
    - \LL(\tilde{\upsilonv}_{\thetavs},\upsilonvs). 
\end{EQA}
Remember that with \(\alp(\upsilonv,\upsilonvs)\) from \eqref{eq: error in taylor approx}
\begin{EQA}[c]
\LL(\upsilonv,\upsilonvs)=\La(\upsilonv,\upsilonvs)+\alp(\upsilonv,\upsilonvs),
\end{EQA}
such that
\begin{EQA}[c]
\LL(\tilde{\upsilonv},\upsilonvs)
    - \LL(\tilde{\upsilonv}_{\thetavs},\upsilonvs)=\La(\tilde{\upsilonv},\upsilonvs)
    - \La(\tilde{\upsilonv}_{\thetavs},\upsilonvs) +\alp(\tilde{\upsilonv},\upsilonvs)
    - \alp(\tilde{\upsilonv}_{\thetavs},\upsilonvs).
\end{EQA}
By the quadratic structure it holds true that
\begin{EQA}
\La(\tilde{\upsilonv},\upsilonvs)&=&\La\left(\upsilonvs+\DF^{-2}\nabla,\upsilonvs \right)+ \|\DF(\tilde{\upsilonv}-\upsilonvs)-\DF^{-1}\nabla\|^{2}/2,\\
\La(\tilde{\upsilonv}_{\thetavs},\upsilonvs)&=&\La\left(\upsilonvs+\HH^{-2}\nabla_{\etav},\upsilonvs \right)+ \|\HH(\tilde \etav_{\thetavs}-\etavs)-\HH^{-1}\nabla_{\etav}\|^{2}/2.
\end{EQA}
Equations \eqref{eq: refined bound for mle} and \eqref{eq: refined bound for thetavs mle} give \( \LCS(\rups,\xx)\subseteq\{ \tilde{\upsilonv}, \tilde{\upsilonv}_{\thetavs} 
\in \Upss(\rupf)\}\cap\bigl\{ \| \DF^{-1} \score \| \le \rupf, \, 
\| \HH^{-1} \score_{\etav} \| \le  \rupf \bigr\}\subset\Omega \). This implies with Theorem~\ref{theo: local expansion full parameter} (adapted appropriately for \(\tilde\etav_{\thetavs}\in\R^{\dimh}\)) that
\begin{EQA}
\|\DF(\tilde{\upsilonv}-\upsilonvs)-\DF^{-1}\nabla\|^{2}\le \Excgr(\rupf,\xx)^{2}, && \|\HH(\tilde \etav_{\thetavs}-\etavs)-\HH^{-1}\nabla_{\etav}\|^{2}\le \Excgr(\rupf,\xx)^{2}.
\end{EQA}
With Lemma~\ref{LLabmProjse} we find
\begin{EQA}
\La\left(\upsilonvs+\DF^{-2}\nabla,\upsilonvs \right)-\La\left(\upsilonvs+\HH^{-2}\nabla_{\etav},\upsilonvs \right)&=& \|\DF^{-1}\nabla\|^{2}/2-\|\HH^{-1}\nabla_{\etav}\|^{2}/2\\	
	&=&\sup_{\upsilonv} \La(\upsilonv,\upsilonvs)
    - \sup_{\upsilonv: \, \Proj\upsilonv=\thetavs} \La(\upsilonv,\upsilonvs)\\
    &=&
    \frac{1}{2} \bigl\| \DPr^{-1} \scorer_{\thetav} \bigr\|^{2}.
\end{EQA}
Consequently we can infer with Lemma~\ref{lem: distance mle and cond mle} and Lemma~\ref{lem: bound on alpha distance} that on \( \LCS(\rups,\xx)\subset\Omega\)
\begin{EQA}
|2\Lr(\tilde{\thetav},\thetavs)-\|\DPr^{-1}\scorer\|^{2}|&=&2\left|\LL(\tilde{\upsilonv},\upsilonvs)
    - \LL(\tilde{\upsilonv}_{\thetavs},\upsilonvs)- \frac{1}{2} \bigl\| \DPr^{-1} \scorer_{\thetav} \bigr\|^{2}\right|\\
    &\le& 2\Excgr(\rupf,\xx)^{2}+2 |\alp(\tilde{\upsilonv},\upsilonvs)
    - \alp(\tilde{\upsilonv}_{\thetavs},\upsilonvs)|\\    	
    &\le&  2(\rr_{\dimp}+\Excgr(\rupf,\xx)) \Excgr(\rupf,\xx).
\end{EQA}
The proof of \eqref{eq: Wilks in main theo} on \( \LCS(\rups,\xx)\subset\Omega \) is completed after substituting \(\rupf=\zz(\xx,\BB)/(1-\corrDF)+\Excgr(\rups,\xx)\) and redefining \(\rr_{\dimp}\in\R\) as in the formulation of the theorem.


\subsubsection{\( \LCS(\rups,\xx)\subset\Omega \) is of dominating probability}
We show that 
the set \( \LCS(\rups,\xx)\subset\Omega \) from \eqref{eq: def of local set} is of dominating probability.
By definition 
\begin{EQA}
    \LCS(\rups,\xx)
    \eqdef& \bigl\{ &
    \tilde{\upsilonv},\tilde{\upsilonv}_{\thetavs}\in\Upss(\rups),
    \\
    &&    \| \DF^{-1} \score \| \le \zz(\xx,\BB),\, 
        \| \HH^{-1} \score_{\etav} \| \le \zz(\xx,\BB)/(1-\corrDF),\, \|\DPr^{-1}\scorer\|\le \zz(\xx,\BBr) 
    \bigr\}\\
    &\cap&\bigg\{\sup_{\upsilonv\in\Upss(\rr)}\|\UU(\upsilonv)\|\le 6\nuno \rhor \zzQ(\xx,\dimtotal) \bigg\}.
\end{EQA}
By the definition of \(\rups>0\) we have
\begin{EQA}[c]
    \P\bigl\{  \tilde{\upsilonv},  \tilde{\upsilonv}_{\thetavs} \not\in \Upss(\rups) \bigr\}
    \le 
    \ex^{-\xx}.
\label{tutusxx}
\end{EQA} 
Further we can use \( (\bb{\AssId})_{\dimh} \) to find
\begin{EQA}[c]
     \| \HH^{-1} \score_{\etav} \|^{2}\le \|\DP^{-1} \score_{\thetav}\|^{2}+\|\HH^{-1} \score_{\etav}\|^{2}\le\frac{1}{1-\corrDF^{2}}\|\DF^{-1}\score\|^{2},
\end{EQA}
which implies that
\begin{EQA}[c]
    \{ \| \DF^{-1} \score \| \le \zz(\xx,\BB)\}\subseteq \{\| \HH^{-1} \score_{\etav} \| \le  \zz(\xx,\BB)/(1-\corrDF)\}.
\end{EQA}  
To control the probability \( \P\bigl( \| \DF^{-1} \score \| > \zz(\xx,\BB) \bigr) \) 
we apply Proposition~\ref{theo: dev bounds quad forms} with 
\begin{EQA}
    \BB
    &=&
    \DF^{-1} \VFc^{2} \DF^{-1}.    
\end{EQA}
With the definitions from Section~\ref{ap: deviation of quadratic forms} we have
\begin{EQA}
    \P\left( \| \DF^{-1} \score \| > \zz(\xx,\BB)\right)    
    & \le &
    2\ex^{-\xx} ,\quad
    \P\bigl( \| \DPr^{-1} \scorer \| > \zz(\xx,\BBr)\bigr)
    \le 
    2\ex^{-\xx} .
\end{EQA}
Further
\begin{EQA}
&&\nquad\P\left(\sup_{\upsilonv\in\Upss(\rr)}\|\UU(\upsilonv)\|\le 6\nuno \rhor \zzQ(\xx,\dimtotal)\right)\ge 1-\ex^{-\xx},
\end{EQA}
by Theorem~\ref{Tsqnorm} which is applicable because \( \bb{(\CS \DF_{1})}\) implies \eqref{eq: zsmu} with \(\|\cdot\|_{\UU}=\|\DF(\cdot)\|\).

Together these bounds yield
\begin{EQA}[c]
    1 - \P\bigl( \LCS(\rups,\xx) \bigr)
    \le 
    6 \ex^{-\xx}.
\label{1PLCSrdrr}
\end{EQA}

\subsection{Proof of Proposition~\ref{lem: large dev refinement}}
Remember that on the set \(\LCS(\rups,\xx)\subset \Omega\) from~\ref{eq: def of local set} we have 
\begin{EQA}[c]
\label{eq: full expansion}
\|\DF(\tilde{\upsilonv}-\upsilonvs)-\DF^{-1}\score\|\le \Excgr(\rups,\xx).
\end{EQA}
With the triangular inequality this gives
\begin{EQA}[c]
\|\DF(\tilde{\upsilonv}-\upsilonvs)\|\le \|\DF^{-1}\score\|+ \Excgr(\rups,\xx).
\end{EQA}
Now on \( \LCS(\rups,\xx) \) we have \( \| \DF^{-1} \score \| \le  \zz(\xx,\BB)\), which implies
\begin{EQA}[c]
\label{eq: refined bound for mle}
\|\DF(\tilde{\upsilonv}-\upsilonvs)\|\le \zz(\xx,\BB)+ \Excgr(\rups,\xx).
\end{EQA}
With the same arguments we find \(\|\DF(\tilde{\upsilonv}_{\thetavs}-\upsilonvs)\|\le \rupf\).
The claim follows as in the proof of Theorem~\ref{theo: main theo finite dim} with \(\rups>0\) replaced with \(\rupf>0\).

%
%
%

\subsection{Proof of Proposition~\ref{lem: critical dim in fisher}}
\begin{proof}
The profile MLE can be calculated easily
\begin{EQA}[c]
\tilde \thetav=f^{-1}\left(\frac{1}{n}\sum_{i=1}^{n}\Yv_i\right)=f^{-1}\left(\frac{1}{n}\sum_{i=1}^{n}\varepsilonv_i\right)
		  =\frac{1}{n}\sum_{i=1}^{n}(\varepsilonv_{i})_1-
		      \|\frac{1}{n}\sum_{i=1}^{n}\varPi_{\etav}\varepsilonv_i\|^{2}.
\end{EQA}
We can show, that the conditions of Section~\ref{sec: conditions} are satisfied and that
\(
\DF^{2}=nI_{\dimn},
\) such that \(\DPr^{2}=n\) and \(\xivr=\frac{1}{\sqrt{\nsize}}\sum_{i=1}^{n}(\varepsilonv_i)_1\). But we immediately see that with a standard normal random variable \(\Zv_{\dimn-1}\in\R^{\dimn-1}\)
\begin{EQA}[c]
\|\DPr(\tilde \thetav-\thetavs)-\xivr\|=\|\sqrt{\nsize}(\tilde \thetav-\thetavs)-\xivr\|=\|\sqrt{\nsize}\tilde \thetav-\xivr\|=\frac{1}{\sqrt{\nsize}}\|\Zv_{\dimn-1}\|^{2}\sim \chi^{2}_{\dimn-1}/\sqrt{\nsize}.
\end{EQA}
This means that if \(\dimn= O(n^{1/2})\) the estimator is not root-n consistent. 
For \(\sqrt{\nsize}=o(\dimn)\) the root-n bias goes to infinity almost surely. 
Clearly if \(\dimn= o(n^{1/2})\) the Fisher expansion is accurate.

Concerning the Wilks phenomenon note that 
\begin{EQA}[c]
\LL(\tilde{\upsilonv})=0.
\end{EQA}
From first order criteria we derive that for all \(j=1,\ldots,\dimh\)
\begin{EQA}
0&=&\bar \varepsilon_j-(\bar\varepsilon_1^{2}-\bar\varepsilon_1\|\tilde\etav_{\thetavs}\|^{2}-1)(\tilde{\eta}_{\thetavs})_j,\\
\|\tilde\etav_{\thetavs}\|^{2}&=&\|\bar\varepsilonv_{\etav}\|^{2}/(\bar\varepsilon_1^{2}-\bar\varepsilon_1\|\tilde\etav_{\thetavs}\|^{2}-1)^{2},
\end{EQA}
where \(\bar\varepsilon_j=\frac{1}{n}\sum_{i=1}^{n}(\varepsilonv_i)_j\) and \(\bar\varepsilon_{\etav}=\frac{1}{n}\sum_{i=1}^{n}\Pi_{\etav}\varepsilonv_i \). As \(\bar\varepsilon_1\to 0\) almost surely we can assume that
\begin{EQA}[c]
\|\tilde\etav_{\thetavs}\|^{2}\cong\|\bar\varepsilonv_{\etav}\|^{2}.
\end{EQA}
This gives
\begin{EQA}
2\Lr(\tilde\thetav,\thetavs)&=&n(\bar\varepsilon_1-\|\tilde\etav_{\thetavs}\|^{2})^{2}+n\|\bar\varepsilonv_{\etav}-\tilde\etav_{\thetavs}\|^{2}\\
	&=&n\bar\varepsilon_1^{2}-2n \bar\varepsilon_1 \|\tilde\etav_{\thetavs}\|^{2}+ n  \|\tilde\etav_{\thetavs}\|^4+n\frac{\bar\varepsilon_1^{2}-\bar\varepsilon_1\|\tilde\etav_{\thetavs}\|^{2}}{\bar\varepsilon_1^{2}-\bar\varepsilon_1\|\tilde\etav_{\thetavs}\|^{2}-1}\|\tilde\etav_{\thetavs}\|^{2}\\
	&\cong&n\bar\varepsilon_1^{2}-2n \bar\varepsilon_1 \|\bar\varepsilonv_{\etav}\|^{2}+ n  \|\bar\varepsilonv_{\etav}\|^4+n\frac{\bar\varepsilon_1^{2}-\bar\varepsilon_1\|\bar\varepsilonv_{\etav}\|^{2}}{\bar\varepsilon_1^{2}-\bar\varepsilon_1\|\bar\varepsilonv_{\etav}\|^{2}-1}\|\bar\varepsilonv_{\etav}\|^{2}\\
	&\cong&n\bar\varepsilon_1^{2}-2n \bar\varepsilon_1 \|\bar\varepsilonv_{\etav}\|^{2}+ n  \|\bar\varepsilonv_{\etav}\|^4.
\end{EQA}
Now if \(\dimn^{2}/n\to 0\) we obtain the Wilks phenomenon. Otherwise if \(\dimn^{2}/n\to\CONST\) the last term on the right hand side is not a chi square distribution with one degree of freedom, while it tends to \(\infty\) almost surely if \(\dimn^{2}/n\to\infty\).
\end{proof}

\subsection{Proof of Theorem~\ref{theo: wilks bias}}

Remember
\begin{EQA}
\tilde{\upsilonv}_{\thetavs_{\dimh},\dimh}\eqdef\left(\thetavs_{\dimh},\ \argmax_{\etav\in\R^{\dimh} } \LL_{\dimh}(\thetavs_{\dimh},\etav)\right),&&
\tilde{\upsilonv}_{\thetavs,\dimh}\eqdef\left(\thetavs,\ \argmax_{\etav\in\R^{\dimh} } \LL_{\dimh}(\thetavs,\etav)\right).
\end{EQA}
Define
\begin{EQA}[c]
	\deltar_{\dimp}^{*}
	\eqdef 
	\left(2+\frac{1+\corrDF}{1-\corrDF}\right)
	\Excgr\bigl(\rr_2,\xx\bigr) + \frac{1+\corrDF}{1-\corrDF} \alpha(\dimh),
\end{EQA}
and for \(\rr_2\eqdef \rups+\deltar_{\dimp}^{*}\)
\begin{EQA}
\mathcal A(\xx,\rr_2)&\eqdef&\bigg\{ 
    \tilde{\upsilonv}, \tilde{\upsilonv}_{\thetavs_{\dimh}}, \tilde{\upsilonv}_{\thetavs}\in\Upss(\rr_2),\,
    \|\DPr^{-1}\scorer\| \le  \zz(\xx,\BBr),\, 
    \| \HH^{-1} \score_{\etav} \| \le  \zz(\xx,\BB)/(1-\corrDF),\\
    &&\| \DF^{-1} \score \| \le \zz(\xx,\BB) 
    \bigg\}\cap\bigg\{\sup_{\upsilonv\in\Upss(\rr)}\|\UU(\upsilonv)\|\le 6\nuno \rhor \zzQ(\xx,\dimtotal) \bigg\}\subset \Omega,
\end{EQA}
with \(\UP(\upsilonv)\in\R^{\dimtotal}\) from \eqref{eq: def of U process}. We use similar arguments as in the proof of Theorem~\ref{theo: main theo finite dim}. First we write
\begin{EQA}
\LL_{\dimh}\Big(\tilde{\upsilonv}_{\thetavs_{\dimh}},\tilde{\upsilonv}_{\thetavs}\Big)=\La\Big(\tilde{\upsilonv}_{\thetavs_{\dimh}},\upsilonvs_{\dimh}\Big)-\La\Big(\tilde{\upsilonv}_{\thetavs},\upsilonvs_{\dimh}\Big)+ \alp\Big(\tilde{\upsilonv}_{\thetavs_{\dimh}},\upsilonvs_{\dimh}\Big)-\alp\Big(\tilde{\upsilonv}_{\thetavs},\upsilonvs_{\dimh}\Big),
\end{EQA}
where
\begin{EQA}
\La\Big(\upsilonv,\upsilonvs_{\dimh}\Big)&=&\nabla_{\dimp+\dimh}\LL(\upsilonvs_{\dimh})(\upsilonv-\upsilonvs_{\dimh})-\|\DF_{\dimh}(\upsilonv-\upsilonvs_{\dimh})\|^{2}/2.
\end{EQA}
We infer by the quadratic structure of \(\La(\cdot)\)
\begin{EQA}
\La\Big(\tilde{\upsilonv}_{\thetavs_{\dimh}},\upsilonvs_{\dimh}\Big)&=& \sup_{\substack{\upsilonv\in\R^{\dimtotal}\\ \upsilonv=(\thetavs_{\dimh},\etav)}}\La\Big(\upsilonv,\upsilonvs_{\dimh}\Big)+ \|\HH_{\dimh}(\tilde{\etav}_{\thetavs_{\dimh}}-\etavs_{\dimh})-\HH_{\dimh}^{-1}\nabla_{\etav}\|^{2}/2,\\
\La\Big(\tilde{\upsilonv}_{\thetavs},\upsilonvs_{\dimh}\Big)&=& \sup_{\substack{\upsilonv\in\R^{\dimtotal}\\ \upsilonv=(\thetavs,\etav)}}\La\Big(\upsilonv,\upsilonvs_{\dimh}\Big)+ \|\HH_{\dimh}(\tilde{\etav}_{\thetavs}-\etavs_{\dimh})-\HH_{\dimh}^{-1}(\nabla_{\etav}-\A_{\dimh}(\thetavs-\thetavs_{\dimh}))\|^{2}/2.
\end{EQA}
As in the proof of Theorem~\ref{theo: main theo finite dim} we find
\begin{EQA}[c]
\label{eq: bound for hh dif of eta mle to eta score}
\|\HH_{\dimh}(\tilde{\etav}_{\thetavs_{\dimh}}-\etavs_{\dimh})-\HH_{\dimh}^{-1}\nabla_{\etav}\|^{2}/2\le \Excgr(\rups,\xx)^{2}/2.
\end{EQA}
Further
\begin{EQA}
&&\nquad\|\HH_{\dimh}(\tilde{\etav}_{\thetavs}-\etavs_{\dimh})-\HH_{\dimh}^{-1}(\nabla_{\etav}-\A_{\dimh}(\thetavs-\thetavs_{\dimh}))\|\\
	&\le& \|\HH_{\dimh}(\tilde{\etav}_{\thetavs}-\tilde{\etav}_{\thetavs_{\dimh}})\|
	+\|\HH_{\dimh}(\tilde{\etav}_{\thetavs_{\dimh}}-\etavs_{\dimh})-\HH_{\dimh}^{-1}(\nabla_{\etav}-\A_{\dimh}(\thetavs-\thetavs_{\dimh}))\|\\
	&&+\corrDF \|\DP(\thetavs_{\dimh}-\thetavs)\|.
\end{EQA}
With equation \eqref{eq: bound for hh dif of eta mle to eta score}, equation \eqref{eq: bound for hh diff eta mle cond and non cond} from the proof of Lemma~\ref{lem: distance mle and cond mle} and with assumption \(\bb{(bias)}\) this gives on \(\mathcal A(\xx,\rr_2)\)
\begin{EQA}[c]
\|\HH_{\dimh}(\tilde{\etav}_{\thetavs}-\etavs_{\dimh})-\HH_{\dimh}^{-1}(\nabla_{\etav}-\A_{\dimh}(\thetavs-\thetavs_{\dimh}))\|\le 3\Excgr\bigl(\rr_2,\xx\bigr)+\frac{2\corrDF}{1-\corrDF}\alpha(\dimh).
\end{EQA}
Again with the arguments in the proof of Lemma~\ref{lem: distance mle and cond mle} we find
\begin{EQA}[c]
\|\DF(\tilde{\upsilonv}_{\thetavs}-\tilde{\upsilonv}_{\thetavs_{\dimh}})\|\le \left(2+\frac{1+\corrDF}{1-\corrDF}\right)\Excgr\bigl(\rr_2,\xx\bigr)+\frac{1+\corrDF}{1-\corrDF} \alpha(\dimh)= \deltar_{\dimp}^{*}.
\end{EQA}
With Lemma~\ref{lem: bound on alpha distance} this implies
\begin{EQA}[c]
	\left|
		\alp\Big(\tilde{\upsilonv}_{\thetavs_{\dimh}},\upsilonvs_{\dimh}\Big)
		- \alp\Big(\tilde{\upsilonv}_{\thetavs},\upsilonvs_{\dimh}\Big)
	\right|
	\le 
	\deltar_{\dimp}^{*}\Excgr\bigl(\rr_2,\xx\bigr).
\end{EQA}
Together this gives
\begin{EQA}
	&&\nquad
	\left |\LL_{\dimh}\Big(\tilde{\upsilonv}_{\thetavs_{\dimh}},\tilde{\upsilonv}_{\thetavs}\Big) 
	-\left(\sup_{\substack{\upsilonv\in\R^{\dimtotal} \\ \upsilonv = (\thetavs_{\dimh},\etav)}}
	\La\Big(\upsilonv,\upsilonvs_{\dimh}\Big) 
	- \sup_{\substack{\upsilonv\in\R^{\dimtotal} \\ \upsilonv=(\thetavs,\etav)}}
		\La\Big(\upsilonv,\upsilonvs_{\dimh}\Big)\right)\right|
	\\
	&\le&  
	\deltar_{\dimp}^{*}\Excgr\bigl(\rr_2,\xx\bigr)+{\deltar_{\dimp}^{*}}^{2}/2
	+ \Excgr\bigl(\rr_2,\xx\bigr)^{2}/2=\left(\deltar_{\dimp}^{*}
	+ \Excgr\bigl(\rr_2,\xx\bigr)\right)^{2}/2.
\end{EQA}
We write with \(\nabla_{\etav}=(\nabla_{\eta_1},\dots,\nabla_{\eta_{\dimh}})\)
\begin{EQA}
\sup_{\substack{\upsilonv\in\R^{\dimtotal}\\\upsilonv=(\thetavs,\etav)}}\La\Big(\upsilonv,\upsilonvs_{\dimh}\Big)&=&\|\HH_{\dimh}^{-1/2}\nabla_{\etav}\LL(\upsilonvs)\|^{2}/2+\|\DPr(\thetavs_{\dimh}-\thetavs)\|^{2}/2\\
	&&+(\thetavs_{\dimh}-\thetavs)^{\T}\left(\nabla_{\thetav}\LL-\A_{\dimh}\HH_{\dimh}^{-1}\nabla_{\etav}\LL \right)\\
	&=&\|\HH_{\dimh}^{-1/2}\nabla_{\etav}\LL(\upsilonvs)\|^{2}/2+r,
\end{EQA}
where by \(\bb{(bias)}\) and Proposition~\ref{theo: dev bounds quad forms}
\begin{EQA}[c]
|r|\le \|\DPr_{\dimh}(\thetavs_{\dimh}-\thetavs)\|^{2}+\|\DPr_{\dimh}(\thetavs_{\dimh}-\thetavs)\|\|\DPr_{\dimh}^{-1}\breve\nabla_{\thetav}\|	\le \alpha^{2}(\dimh)+\alpha(\dimh)\zz(\xx,\BBr).
\end{EQA}
Further
\begin{EQA}[c]
\sup_{\substack{\upsilonv\in\R^{\dimtotal}\\\upsilonv=(\thetavs_{\dimh},\etav)}}\La\Big(\upsilonv,\upsilonvs_{\dimh}\Big)=\|\HH_{\dimh}^{-1/2}\nabla_{\etav}\LL(\upsilonvs)\|^{2}/2.
\end{EQA}
Consequently on the set \(\mathcal A(\xx)\subset \Omega\)
\begin{EQA}
&&\nquad\bigl| 2 \max_{\etav\in \R^{\dimh}}\LL_{\dimh}(\thetavs_{\dimh},\etav)- \max_{\etav\in \R^{\dimh}}2 \LL_{\dimh}(\thetavs,\etav)\bigr|= \left|\LL_{\dimh}\Big(\tilde{\upsilonv}_{\thetavs_{\dimh}},\tilde{\upsilonv}_{\thetavs}\Big)\right|\\
	&\le& \left(\deltar_{\dimp}^{*}+\Excgr\bigl(\rr_2,\xx\bigr)\right)^{2}+2\alpha^{2}(\dimh)+2\alpha(\dimh)\zz(\xx,\BBr).
\end{EQA}
The claim follows because the result \eqref{eq: Wilks in main theo} of Theorem~\ref{theo: main theo finite dim} occurs on \(\mathcal A(\xx,\rr_2)\subset\Omega\). It remains to note that the set \(\mathcal A(\xx,\rr_2)\subset \Omega\) is of probability greater  \(1-6\ex^{-\xx}\) by the choice of \(\rr_2>0\). We see this with the same arguments as in the proof of Lemma~\ref{lem: distance mle and cond mle}
\begin{EQA}
&&\nquad\P\{\tilde{\upsilonv}, \tilde{\upsilonv}_{\thetavs_{\dimh}},\tilde{\upsilonv}_{\thetavs}\in\Ups_{0,\dimh}(\rr_2)\cap \LCS(\rups,\xx)\}\\
&\ge& \P\left\{\tilde{\upsilonv}, \tilde{\upsilonv}_{\thetavs_{\dimh}}\in\Ups_{0,\dimh}\left(\rr_2-\|\DF(\tilde{\upsilonv}_{\thetavs}-\tilde{\upsilonv}_{\thetavs})\|\right)\cap\LCS(\rups,\xx)\right\}\\
	&\ge& \P\left\{\tilde{\upsilonv},\tilde{\upsilonv}_{\thetavs_{\dimh}}\in\Ups_{0,\dimh}\left(\rr_2-\deltar_{\dimp}^{*}\right)\cap\LCS(\rups,\xx)\right\}\\
	&=&\P\{\LCS(\rups,\xx)\}\ge 1-6\ex^{-\xx}.
\end{EQA}

\subsection{Proof of Corollary~\ref{cor: efficiency in iid case} and \ref{cor: wilks in iid case}}

We will only prove Corollary~\ref{cor: efficiency in iid case}, as the the proof of Corollary~\ref{cor: wilks in iid case} is very similar.

Define
\begin{EQA}[c]
\VF_{\dimh}^{2}=\Cov\bigl(\score_{\dimp+\dimh} \LL_{\dimh}(\upsilonvs_{\dimh})\bigr),\,{\BB}_{\dimh}=\DF_{\dimh}^{-1}\VF_{\dimh}^{2}\DF_{\dimh}^{-1},\\
\scorer_{\thetav,\dimh}
    =
    \score_{\thetav} - \A_{\dimh} \HH_{\dimh}^{-2} \score_{\etav},\, \VPr_{\dimh}^{2}=\Cov(\scorer_{\thetav}),\, \breve{\BB}_{\dimh}=\DPr_{\dimh}^{-1}\VPr^{2}_{\dimh}\DPr_{\dimh}^{-1}.
\end{EQA}

Remember \( \dimtotal = \dimp + \dimh\in\N \) and that the point \( \upsilonvs\in l^{2} \) is defined by maximizing the expected 
log-likelihood for the sieved functional models \(\LL_{\dimh}\) and the operators \(\DF_{\dimh}^{2}\in\R^{\dimtotal\times\dimtotal}, \DF^{2}\in L(l^{2},l^{2}) \) correspond to this point.

Fix \(\xx \) and define \( \rups \) by 
\(\P\big\{\tilde{\upsilonv}\in\Upss(\rr_{0,n})\big\}\ge 1-\ex^{\xx}\to 1 \). 
Then we get with Theorem ~\ref{theo: main theo finite dim} applied to 
\( \tilde{\thetav}_{\dimh} \) from \eqref{ttdhttdhs} that with probability greater \(1-6 \ex^{-\xx}\)
\begin{EQA}[c]
\label{eq: convergence of estimator for dimh}
    \|\DPrp \bigl( \tilde{\thetav}_{\dimh} - \thetavs_{\dimh} \bigr)
    - \xivr_{\dimh}(\upsilonvs_{\dimh})\|\le \Excgr(\rups,\xx).
\label{DPrctthsh}
\end{EQA}  
We write
\begin{EQA}
	&& \nquad
	\DPr \bigl( \tilde{\thetav}_{\dimh} - \thetavs \bigr)    - \xivr_{\dimh}(\upsilonvs_{\dimh})
	\\
	&=&
	\DPrp \bigl( \tilde{\thetav}_{\dimh} - \thetavs_{\dimh} \bigr)- \xivr_{\dimh}(\upsilonvs_{\dimh})
	+ ( \DPrp-\DPr)\bigl( \tilde{\thetav}_{\dimh} - \thetavs_{\dimh} \bigr) 
	+ \DPrp \bigl( \thetavs_{\dimh} - \thetavs \bigr).
\end{EQA}
By \eqref{eq: convergence of estimator for dimh} it suffices to bound
\(\| (\DPrp-\DPr)( \tilde{\thetav}_{\dimh} - \thetavs_{\dimh} )\|\) and \( \|\DPrp(\thetavs_{\dimh} - \thetavs) \|\). With assumption \(\bb{(bias)}\) and \(\bb{(bias')}\). We get
\begin{EQA}
\|\DPr \bigl( \tilde{\thetav}_{\dimh} - \thetavs \bigr)    - \xivr_{\dimh}(\upsilonvs_{\dimh})\|&\le &\Excgr(\rups,\xx)+o(1)\|\DPrp \bigl(\tilde{\thetav}_{\dimh} - \thetavs \bigr)\|+\alpha(\dimh),
\end{EQA}
where we used that
\begin{EQA}
	&&\nquad
	\|( \DPrp-\DPr)\bigl( \tilde{\thetav}_{\dimh} - \thetavs_{\dimh} \bigr)\|
	\\
	&\le& 
	\|( \DPrp-\DPrp(\upsilonvs))\bigl( \tilde{\thetav}_{\dimh} - \thetavs_{\dimh} \bigr)\|
	+ \|( \DPrp(\upsilonvs)-\DPr)\bigl( \tilde{\thetav}_{\dimh} - \thetavs_{\dimh} \bigr)\|
	\\
	&\le& 
	\|\DPrp \bigl(\tilde{\thetav}_{\dimh} - \thetavs \bigr)\|
	\,
	\Bigl( \| \Id - \DPrp^{-1}\DPrp^{2}(\upsilonvs)\DPrp^{-1}\|^{1/2}
	\\
	&& 
	\quad	+ \, \| \Id - \DPrp(\upsilonvs)^{-1}\DPr^{2}(\upsilonvs)\DPrp(\upsilonvs)^{-1}\|^{1/2}
	\|\DPrp(\upsilonvs)\DPrp^{-1}\|
	\Bigr) .
\end{EQA}
Now we use again \eqref{eq: convergence of estimator for dimh} and the fact, that on the set on which \eqref{eq: convergence of estimator for dimh} holds we have \(\|\xivr\|\le \zz(\xx,\BBr_{\dimh})\) (see \eqref{eq: def of local set}) to obtain 
\begin{EQA}[c]
	\| \DPr_{\dimh} (\tilde{\thetav}_{\dimh} - \thetavs_{\dimh}) \|
	\le 
	\|\xivr\|+\Excgr(\rups,\xx)
	\le 
	\zz(\xx,\BBr_{\dimh})+\Excgr(\rups,\xx).
\end{EQA}
Combining these bounds gives 
\begin{EQA}
	\|\DPr \bigl( \tilde{\thetav}_{\dimh} - \thetavs \bigr)    - \xivr_{\dimh}(\upsilonvs_{\dimh})\|
	&\le &
	\Excgr(\rups,\xx)
	+ o(1)\left(\zz(\xx,\BBr_{\dimh})+\Excgr(\rups,\xx) \right)+\alpha(\dimh).
\end{EQA}

Fix a sequences \( (\dimh) = (\dimh(\nsize)) \), 
\( (\xxn)  \) from \eqref{eq: conditions on approx error in semiefficient} further ensure that \(\xxn\to \infty\) slow enough to ensure that 
\( o(1)\left(\zz(\xx,\BBr_{\dimh})+\Excgr(\rups,\xx) \right)=o(1)\). 
Due to \eqref{eq: conditions on approx error in semiefficient} we have that for any \(\eps>0\) there exists an \(n\in\N\) such that
\begin{EQA}[c]
	\P(\|\DPr(\tilde \thetav_{\dimh}-\thetavs)-\xivr_{\dimh}(\upsilonvs_{\dimh})\|\ge \eps)
	\le 
	1-6\ex^{\xxn}.
\end{EQA}
As \(\xxn\to \infty\) we get the claim by Slutsky's Lemma once we showed that \( \xivr_{\dimh}(\upsilonvs_{\dimh}) \) is asymptotically standard normal.

For this observe
\begin{EQA}
    \xivr_{\dimh}(\upsilonvs_{\dimh})
    &=&
    \DPr_{\dimh}^{-1} (\score_{\thetav} - \A_{\dimh}\HH_{\dimh}^{-2} \score_{\etav}) 
    	\LL(\upsilonvs_{\dimh})
    \\
    &=&
    \frac{1}{\sqrt{\nsize}}\sum_{i=1}^{\nsize} 
        (\frac{1}{\sqrt{\nsize}} \DPr_{\dimh})^{-1}(\score_{\thetav} \lkh_{i}(\upsilonvs_{\dimh}) - \A_{\dimh}\HH_{\dimh}^{-2}\score_{\etav} \lkh_{i}(\upsilonvs_{\dimh}))\\
     &\eqdef&\frac{1}{\sqrt{\nsize}}\sum_{i=1}^{\nsize} \Xv_{i}.
    \end{EQA}
Due to assumptions \(\bb{(bias'')}\) we have \(\Cov(\Xv_i)\to \Id_{\dimp}\) as \( \Cov(\hat \Xv_{i})= \Id_{\dimp}\). Consequently
\begin{EQA}
    \xivr_{\dimh}(\upsilonvs_{\dimh})
    &=&
    \frac{1}{\sqrt{\nsize}}\sum_{i=1}^{\nsize} \Xv_{i},
\end{EQA}
where the random vectors 
\( \Xv_{i} \)
are i.i.d. with zero mean and covariance tending to \( \Id_{\dimp} \), such that by a slightly generalized central limit theorem
\begin{EQA}
\xivr_{\dimh}(\upsilonvs_{\dimh})& \tow &\ND(0,\Id_{\dimp}).
\end{EQA}

\section{A bound for the norm of a random process}
\label{ap: bound for norm of gradient}
We want to derive for a random process \(\UU(\ups)\in\R^{\dimtotal}\) a bound of the kind
\begin{EQA}[c]
	\P\left(\sup_{\ups \in \Upss(\rr)} \| \UU(\ups) \|
	\ge
	\CONST \zzQ(\xx,\dimtotal)\rr\right)
	\le \ex^{-\xx}.
\end{EQA}
In the following we elaborate how to extend the results of the supplement of \cite{SP2011} on empirical processes to this situation without substantial changes to the bounds.

For this let \( \UU(\ups) \) be a smooth centered random vector process with values in 
\( \R^{\dimtotal} \) and let \(\DF: \R^{\dimtotal}\to \R^{\dimtotal}\) be some linear operator. 
We aim at bounding the maximum of the norm \( \| \UU(\ups) \| \)
over a vicinity \( \Upss \) of \( \upss \).
Suppose that \( \UU(\ups) \) satisfies for each \( 0<\rr<\rr^{*} \) and for all pairs \( \ups,\upsd \in \Upss(\rr) = \bigl\{ \ups \in \Ups \colon \| \ups - \upss \| \leq \rr \bigr\}\subset \R^{\dimtotal} \)

\begin{EQA}
\label{eq: zsmu}
    \sup_{\|\uv\|\le 1} \log \E \exp \biggl\{
        \lambda \frac{\uv^{\T}\bigl(\UU(\ups)-\UU(\upsd)\bigr)}{\rhor\|\ups - \upsd \|}
    \biggr\}
    & \le &
    \frac{\nunu^{2} \lambda^{2}}{2} .
\end{EQA}

\begin{remark}
In the setting of Theorem~\ref{theo: main theo finite dim} we have
\begin{EQA}[c]
\UU(\upsilonv)=\DF^{-1}\Big(\nabla\zetav(\upsilonv) - \nabla\zetav(\upsilonvs)\Big),
\end{EQA}
and condition \eqref{eq: zsmu} becomes \(\bb{(\CS \DF_{1})}\) from \ref{sec: conditions}. 
\end{remark}

\begin{theorem}
\label{Tsqnorm}
Let a random \( \dimtotal \)-vector process \( \UU(\ups) \) 
fulfill \( \UU(\upss) = 0 \), \( \E \UU(\ups) \equiv 0 \), 
and the condition \eqref{eq: zsmu} be satisfied. 
Then for each \( \rr > 0 \), on a set of probability greater \(1-\ex^{-\xx}\)
\begin{EQA}
	\sup_{\ups \in \Upss(\rr)} \| \UU(\ups) \|
	& \leq &
	6 \rhor\nuno \zzQ(\xx,\dimtotal)\rr,
\label{bouuvupsrupsdx22}
\end{EQA}
where with \( \gmd = \nunu \gmb \) and for some \(\entrl>0\)
\begin{EQA}[c]
    \zzq(\xx,\entrl)
    \eqdef
    \begin{cases}
    \sqrt{2 (\xx + \entrl)} 
        & \text{if } \sqrt{2 (\xx + \entrl)} \le \gmd, \\
    \gmd^{-1} (\xx + \entrl) + \gmd/2 
        & \text{otherwise} .
  \end{cases}
\label{PUPxxllc}
\end{EQA}   
\end{theorem}

\begin{remark}
Note that the entropy of the original set is multiplied by 2 as the supremum is taken over \(\Upss(\rr)\times B_{\rr}(0)\subset \R^{2\dimtotal}\). So in order to control the norm \( \| \UU(\ups) \|\) one only pays with this factor.
\end{remark}

\begin{proof}
In what follows, we use the representation
\begin{EQA}[c]
    \| \UU(\ups) \|
    =
    \sup_{\|\uv\|\le \rr} \frac{1}{\rr} \uv^{\T} \UU(\ups) .
\label{UU2ups}
\end{EQA}    
This implies
\begin{EQA}[c]
    \sup_{\ups \in \Upss(\rr)} \| \UU(\ups) \|
    =
    \sup_{\ups \in \Upss(\rr)} \sup_{\|\uv\|\le \rr} 
        \frac{1}{\rr}  \uv^{\T} \UU(\ups) .
\label{UU2vups}
\end{EQA}   

Due to Lemma~\ref{lem: exp moments for scalarproduct with gradient} the process \(\UP(\ups,\uv)\eqdef\frac{1}{\rr}  \uv^{\T} \UU(\ups) \) satisfies condition \eqref{lexpuvUps} as process on \(\R^{2\dimtotal}\). This allows to apply Corollary 2.2 of the supplement of \cite{SP2011} to obtain the desired result.
We get on a set of probability greater \(1-\ex^{-\xx}\) 
\begin{EQA}
	\sup_{\ups \in \Upss(\rr)} \,\, \sup_{\|\uv\|\le \rr} \,\, 
		\left\{ 
			\frac{1}{6 \nuno \rr} \uv^{\T} \UU(\ups) 
			\right\}
	& \leq &
	\zz\Big(\xx,\entrl\bigl(\Upss(\rr)\times \{\B_{\rr}(0)\cap\Upss(\rr)\} \bigr)\Big).
\label{bouuvupsdx}
\end{EQA}
The constant \(\entrl\bigl(\Upss(\rr)\times \{\B_{\rr}(0)\cap\Upss(\rr)\}\bigr)>0\) quantifies the complexity of the set \(\Upss(\rr)\times \{\B_{\rr}(0)\cap\Upss(2\rr)\}\subset \R^{\dimtotal}\times \R^{\dimtotal}\). We point out that for compact \(M\subset \R^{\dimtotal}\) we have \(\entrl(M)= 2 \dimtotal\) (see Supplement of \cite{SP2011}, Lemma 2.10). This gives \(\entrl\bigl(\Upss(\rr)\times \{\B_{\rr}(0)\cap\Upss(\rr)\}\bigr)=4\dimtotal\).
\end{proof}

\begin{lemma}
\label{lem: exp moments for scalarproduct with gradient}
 Suppose that \( \UU(\ups) \) satisfies for each \( \ups \in \Upss \) and each 
\( \gammav\in \R^{\dimtotal} \) with 
\( \| \uv \| \le 1 \) and some norm \(\|\cdot\|_{\UU}\)
\begin{EQA}[c]
    \sup_{\ups \in \Upss} \log \E \exp\Bigl\{ 
		\frac{\lambda( \UU(\ups)-\UU(\upsd))^{\T} \uv}{\rhor\|\ups-\upsd\|_{\UU}} 
	\Bigr\} 
    \le 
    \frac{\nunu^{2} \lambda^{2}}{2} , 
    \qquad 
    |\lambda| \leq \gm.
\label{zsmu}
\end{EQA}
Then for any \(\uv_{1},\uv_{2}\in\R^{\dimtotal} \) with \(\|\uv_i\|_{\UU}\le 2\rr\) and \(\|\uv_i\|\le \rr\)
\begin{EQA}
	\log \E \exp\biggl\{ 
		\frac{\lambda}{2 \rr} 
        \frac{(\UU(\ups)^{\T}\uv_{1}-\UU(\upsd)^{\T}\uv_{2})}{\rhor\sqrt{\|\ups-\upsd\|_{\UU}^{2}+\|\uv_{1}-\uv_{2}\|^{2}}} 
	\biggr\}
	& \leq & \frac{\nunu^{2} \lambda^{2}}{2} ,
	\qquad
	\lambda \leq \gm .
\label{lexpuvUps}
\end{EQA}
\end{lemma}
\begin{proof}
By the H\"older inequality and \eqref{zsmu}, it holds for \( \ups,\upsd\in\Upss( \rr) \)
\begin{EQA}
	&& \nquad
	\log \E \exp\biggl\{ 
		\frac{\lambda}{2 \rr} 
         \frac{\uv_1^{\T} \UU(\ups)- \uv_2^{\T} \UU(\upsd)}{\rhor\sqrt{\|\ups-\upsd\|^{2}_{\UU}+\|\uv_1-\uv_2\|^{2}}}
	\biggr\}
	\\
	&= &
	\log \E \exp\biggl\{ 
		\frac{\lambda}{2 \rr} 
         \frac{\uv_1^{\T} \bigl(\UU(\ups)-\UU(\upsd)\bigr)+(\uv_1^{\T}-\uv_2^{\T}) \UU(\upsd)}{\rhor\sqrt{\|\ups-\upsd\|^{2}_{\UU}+\|\uv_1-\uv_2\|^{2}}}
	\biggr\}
	\\
	& \leq &
	\frac{1}{2}\log \E \exp\biggl\{ 
		\frac{\lambda}{ \rr} 
         \frac{\uv_1^{\T} \bigl(\UU(\ups)-\UU(\upsd)\bigr)}{\rhor\|\ups-\upsd\|_{\UU}}
	\biggr\}+\frac{1}{2}\log \E \exp\biggl\{ 
		\frac{\lambda}{ \rr} 
         \frac{(\uv_1^{\T}-\uv_2^{\T}) \UU(\upsd)}{\rhor\|\uv_1-\uv_2\|}
	\biggr\}
	\\
	& \leq &
	\sup_ {\|\uv\|\le 1} \frac{1}{2}\log \E \exp\biggl\{ 
		\lambda 
         \frac{\uv^{\T}\big (\UU(\ups)-\UU(\upsd)\bigr)}{\rhor\|\ups-\upsd\|_{\UU}}
	\biggr\}\\
	 &&+\sup_{\|\uv\|\le 1}\frac{1}{2}\log \E \exp\biggl\{ 
		\lambda 
         \frac{\uv^{\T}\bigl( \UU(\upsd)-\UU(\upss)\bigr)}{\rhor\|\upsd-\upss\|_{\UU}}
	\biggr\}
	\\
	& \leq &
	\frac{\nunu^{2} \lambda^{2}}{2} ,
	\qquad
	\lambda \leq \gm .
\end{EQA}
\end{proof}

\bibliography{semiquellen,exp_ts,listpubm-with-url}
\end{document}